\documentclass{article}
\usepackage{graphicx}
\usepackage{authblk}
\usepackage{hyperref}
\hypersetup{
	pdftitle={Numerical solution of an optimal control problem with probabilistic and almost sure state constraints},
	pdfauthor={C. Geiersbach and R. Henrion and P. Pérez-Aros}
}
\usepackage{amsmath, color, verbatim}
\usepackage{algorithm} 
\usepackage{algorithmic}
\usepackage[numbers,sort]{natbib}
\usepackage{amsfonts, mathtools, enumitem} 
\usepackage{hyperref}
\usepackage{todonotes}
\usepackage{enumitem}
\usepackage{titlesec}
\usepackage{cleveref}
\usepackage[left=2cm,right=2cm,top=2cm,bottom=2cm]{geometry}
\hypersetup{
  colorlinks=true,
}
\usepackage{amssymb}
\usepackage{amsthm}

\newtheorem{theorem}{Theorem}[section]

\newtheorem{lemma}[theorem]{Lemma}
\newtheorem{proposition}[theorem]{Proposition}

\newtheorem{corollary}[theorem]{Corollary}
\newtheorem{remark}[theorem]{Remark}
\newcounter{assumption}
\newtheorem{assumption}[theorem]{Assumption}

\newcommand{\E}{\mathbb{E}}
\newcommand{\pP}{\mathbb{P}}
\newcommand{\R}{\mathbb{R}}

\newcommand{\D}{\text{ d}}

\begin{document}

\title{Numerical solution of an optimal control problem with probabilistic and almost sure state constraints}


\author{Caroline Geiersbach\thanks{Weierstrass Institute, 10117 Berlin, Germany 
  (\texttt{caroline.geiersbach@wias-berlin.de})}
\and Ren\'e Henrion\thanks{Weierstrass Institute, 10117 Berlin, Germany 
  (\texttt{rene.henrion@wias-berlin.de})}
\and Pedro Pérez-Aros \thanks{Instituto de Ciencias de la Ingenier{\'i}a\\ Universidad de O’Higgins, Chile (\texttt{pedro.perez@uoh.cl})
  }}

\date{\today}

\maketitle

\begin{abstract}
We consider the optimal control of a PDE with random source term subject to probabilistic or almost sure state constraints. In the main theoretical result, we provide an exact formula for the Clarke subdifferential of the probability function without a restrictive assumption made in an earlier paper. The focus of the paper is on numerical solution algorithms. As for probabilistic constraints, we apply the method of spherical radial decomposition. Almost sure constraints are dealt with a Moreau--Yosida smoothing of the constraint function accompanied by Monte Carlo sampling of the given distribution or its support or even just the boundary of its support. Moreover, one can understand the almost sure constraint as a probabilistic constraint with safety level one which offers yet another perspective. Finally, robust optimization can be applied efficiently when the support is sufficiently simple. A comparative study of these five different methodologies is carried out and illustrated. 
\end{abstract}

\section{Introduction}
Many physics-based systems that can be described mathematically as optimization problems contain inputs or parameters that are unknown. Ignoring an available model for the uncertain values, for example a probability distribution, can result in severely sub-optimal solutions. In the context of PDE-constrained optimization under uncertainty,  the framework of stochastic optimization has proved useful due to its rich theory and wealth of numerical methods \cite{Conti2009}. The simplest ansatz is to optimize with respect to a desired \textit{average} outcome \cite{Conti2009,Gahururu2022,Geiersbach2021c,Milz2023} of the underlying random PDE, but robust \cite{Kolvenbach2018,Alla2019} and risk-averse formulations have also been proposed for engineering applications \cite{conti2011,Rockafellar2015} with theoretical investigations in \cite{Kouri2016,Kouri2018,Kouri2019a,alphonse2022risk,Geiersbach2022a}. 

In certain applications, additional constraints on the solution to the random PDE are also desirable, leading to \textit{state constraints}. Systems involving an almost sure or robust model for state constraints have recently been investigated in \cite{Geiersbach2021c,Gahururu2022,Geiersbach2022a,geiersbach2023optimality}. Probabilistic constraints offer the possibility to deal with uncertain restrictions in a robust way that is interpretable with respect to probability. They
have been introduced by Charnes et al.~\cite{charnes1958}. Fundamental algorithmic and theoretical contributions are due to Prékopa, see \cite{prekopa1995}. More recent presentations are provided in \cite{Shapiro2009} or \cite{vanAckooij2020}, respectively. In the last years, one may observe a growing interest in probabilistic constraints as part of PDE-constrained optimization, see e.g., \cite{caillau2018,Farshbaf-Shaker2020,Farshbaf-Shaker2018,geletu2020,goettlich2021,Kouri2023,Perez2022,Teka23,Schuster2022}. These works include proposals for numerical approaches as well as structural investigations.

A central challenge in optimization under uncertainty is in the numerical solution. Classical approaches involve discretization of the stochastic space \cite{Benner2016b,Garreis2017}, stochastic collocation \cite{Kouri2013a,Chen2014,Chen2016,Zahr2019}, or using a sample average approximation \cite{Haber2012a,Ali2017a,Guth2021,VanBarel2019a,Milz2023}. Stochastic approximation, which dynamically samples over the course of optimization, has also gained attention \cite{Geiersbach2019a,Martin2021,Martin2021a,Geiersbach2021c,Geiersbach2023,Kouri2023}. Other innovations include the use of surrogate functions constructed using Taylor approximations of the objective and constraint function \cite{Alexanderian2017a,Chen2021a}. Very few approaches exist that can handle state constraints. A Monte Carlo approximation was used in combination with a Moreau--Yosida penalty in \cite{Gahururu2022}. In~\cite{Kouri2023}, almost sure state constraints were relaxed to an expectation constraint and a stochastic approximation approach was proposed. In \cite{antil2023state}, random fields are approximated by the tensor-train decomposition and state constraints are handled using a Moreau--Yosida-type penalty with a softplus approximation for the positive part function. 

The current work is a follow up paper to \cite{geiersbach2023optimality}, where optimality conditions on a risk-averse PDE-constrained optimization problem with uncertain state constraint was considered. Risk aversion was modeled by probabilistic and almost sure constraints. In the current paper, we pick up the same optimization problem (subject to Poisson's equation with a distributed control and an additional random source term on the right-hand side), but focus on its numerical solution. The paper is organized as follows: after presenting the model in Section \ref{model}, we analyze the problem in Section \ref{sec:analysis}. The main result is an improvement of an exact formula for the Clarke subdifferential of the probability function provided in \cite{geiersbach2023optimality} in that it omits a restrictive assumption on boundedness of the set of feasible realizations of the random vector. Section \ref{numprobcons} is then devoted to the numerical solution of the control problem under probabilistic (uniform) state constraints in 1D and 2D. Our approach for dealing with probabilistic constraints is based on the well-studied spherical radial decomposition. In Section \ref{sec:almost-sure}, we pass to almost sure constraints. For their numerical treatment, we follow a different methodology than for probabilistic constraints, namely sampling of the distribution (or its support) and applying a Moreau-Yosida approximation. Four different approaches are compared with a reference solution obtained by robust optimization.
\section{The model}\label{model}
In this paper, we discuss numerical approaches for solving the following risk-averse PDE-constrained optimization problem under uncertainty:
\begin{subequations}
    \begin{alignat}{3}
    \min_{u \in L^2(D)}\, F(u) & & &&   \label{eq:probuniform-problem-a}\\
    \text{s.t.}  \quad -\Delta  y(x,\omega)  &=  u(x) + f(x,\xi (\omega)), & &&\quad x \in D \quad \text{ $\mathbb{P}$-a.s.},  \label{eq:probuniform-problem-b} \\
 y(x,\omega) &=0, & &&\quad x \in \partial D  \quad \text{$\mathbb{P}$-a.s.},  \label{eq:probuniform-problem-c}\\
 \mathbb{P}(y(x,\omega)&\leq \alpha\quad\forall x\in D)\geq p. & &&  \label{eq:probuniform-problem-d}
    \end{alignat}
    \label{eq:probuniform-problem}
\end{subequations}
Here, $D\subseteq\mathbb{R}^d$ ($d=1,2,3$) is an open and bounded set. Moreover, $F\colon L^2(D)\to\mathbb{R}$ is a convex, Fréchet differentiable cost function, $\xi\sim\mathcal{N}(0,\Sigma)$ is a centered $m$-dimensional Gaussian random vector defined on the probability space $(\Omega,\mathcal{F},\mathbb{P})$, $p\in (0,1]$ is some given probability level, and $\alpha\in\mathbb{R}$ is some upper threshold for the random state $y(\cdot,\omega)$. The function $f\colon \mathbb{R}^d\times\mathbb{R}^m\to\mathbb{R}$ is a random source term.
Inequality \eqref{eq:probuniform-problem-d} is a probabilistic constraint expressing the condition that the random state $y$ stays below $\alpha$ uniformly on $D$ with probability at least $p$. Throughout this paper, we will make the following assumptions on the PDE \eqref{eq:probuniform-problem-b}--\eqref{eq:probuniform-problem-c}. In the following, the notation $\text{meas}(\cdot)$ refers to the $d$-dimensional Lebesgue measure.
\begin{assumption}
\label{ass:PDE-standing}
The open and bounded set $D\subseteq\mathbb{R}^d$ ($d=1,2,3$) is of class $S$ \footnote{For the domain, we use the terminology from \cite{Kinderlehrer1980} and note that class $S$ covers many cases, including Lipschitz or convex domains.}, meaning that there exist constants $\gamma \in (0,1)$ and $r_0 >0$ such that $\text{meas}(B_{r}(x)\backslash D) \geq \gamma \text{meas}( B_r(x))$ for all $x \in \partial D$ and for all $ r< r_0$. Additionally, $u \in L^2(D)$ and the function $f\colon \mathbb{R}^d\times\mathbb{R}^m\to\mathbb{R}$ is defined by 
\begin{equation}
\label{eq:f}
 f(x,z):= f_0(x) + \sum_{i=1}^m z_i\phi_i(x)
\end{equation}
for some given $f_0, \phi_i\in L^2(D)$. 
\end{assumption}
Note that $\mathbb{E}f(\cdot,\xi )=f_0$ on account of $\xi$ being centered. We emphasize that the basic structure of the optimization problem introduced above is intentionally kept simple in order to focus the attention on the aspect of probabilistic state constraints. It is, however, no problem to pass to more general settings, such as: alternative multivariate distributions (Gaussian-like, elliptically symmetric), additional simple control constraints, two-sided state constraints possibly with functional threshold $\alpha$, etc. We will occasionally pick up some of these aspects in subsequent sections.

\section{Analytical properties of the problem}
\label{sec:analysis}
\subsection{General statements on optimization problems with probabilistic constraints}
We start by embedding problem \eqref{eq:probuniform-problem} into a more general framework, which is given by
\begin{equation}\label{optprob}
\min\limits_{u\in U}\,\,F(u)\mbox{ s.t. }\varphi(u)\geq p\quad (p\in (0,1]).
\end{equation}
Here, $U$ is a reflexive and separable Banach space, $F\colon U\to\mathbb{R}$ is some convex, Fréchet differentiable cost function, and $\varphi:U\to\mathbb{R}$ denotes a probability function defined by 
    $\varphi (u):=\mathbb{P}(\omega\mid g(u,\xi(\omega))\leq 0)$. In this last expression, $\xi$ is an $m$-dimensional Gaussian random vector defined on a probability space $(\Omega, \mathcal{F}, \pP)$ and having a centered Gaussian distribution $\mathcal{N}(0,\Sigma)$
with covariance matrix $\Sigma$,
and $g\colon U\times\mathbb{R}^m\to\mathbb{R}$ is some constraint function. 
We make the following general assumption on $g$:
\begin{equation}\tag{GA}
g\mbox{ is locally Lipschitzian and } g(u,\cdot )\mbox{ is convex for all } u\in U. 
\end{equation}
We shall say that $g$ satisfies the {\it condition of moderate growth} at $\bar{u}\in U$, if 
\begin{eqnarray}
&\exists l>0\,\,\forall d\in U:g^{\circ}(\cdot ,z)(u;d)\leq l \left\Vert z\right\Vert^{-m}\exp\left(\frac{\left\Vert z\right\Vert^2}{2\left\Vert \Sigma^{1/2}\right\Vert^2} \right)\left\Vert d\right\Vert&\label{growth}\\ &\forall
u\in \mathbb{B}_{1/l}\left(\bar{u} \right), \forall z: \left\Vert z\right\Vert\geq l&\notag
\end{eqnarray}

Here, $g^{\circ}(\cdot ,z)(u;d)$ refers to the Clarke directional derivative of the locally Lipschitzian (by (GA)) partial function $g(\cdot,z)$ at the argument $u$ in direction $d$. Moreover, $\Sigma^{1/2}$ denotes a root of $\Sigma$. As a consequence of the {\em spherical radial decomposition} of Gaussian random vectors, the total probability function $\varphi$ can be represented as a spherical integral with respect to the uniform distribution $\mu_\zeta$ on $\mathbb{S}^{m-1}$
\begin{equation}\label{srd}
\varphi (u)=\int\limits_{\mathbb{S}^{m-1}}e(u,v) \D \mu_\zeta (v)\quad (u\in U)
\end{equation}
over a one-dimensional radial probability function $e\colon U\times\mathbb{S}^{m-1}\to\mathbb{R}$ defined by 
\begin{equation}\label{radprob}
e(u,v):=\mu_\chi (\{r\geq 0\mid g(u,r\Sigma^{1/2}v)\leq 0\}),
\end{equation}
where $\mu_\chi$ is the one-dimensional Chi- distribution with $m$ degrees of freedom. We refer to \cite{vanAckooij_Henrion_2014,vanAckooij2020b, vanAckooij2020,Hantoute2019,MR4000225,MR4382608} for more details and generalizations of such decomposition to other classes of distributions. 
The following result on subdifferentiation under the integral sign holds true:
\begin{theorem}[\cite{Hantoute2019}, Theorem 5, Corollary 2 and Proposition 6]\label{probderiv}
Under the basic assumptions (GA), let $\bar{u}\in U$ be given such that $g(\bar{u},0)<0$ and \eqref{growth} is satisfied. Then, $\varphi$ and the $e(\cdot ,v)$, $(v\in\mathbb{S}^{m-1})$, are locally Lipschitzian around $\bar{u}$ and, with $\partial^C$ denoting the Clarke subdifferential, it holds that
\begin{equation}\label{clarkeinclu}
\partial^C\varphi (\bar{u})\subseteq\int\limits_{\mathbb{S}^{m-1}}\partial^C_ue(\bar{u},v) \D \mu_\zeta (v).
\end{equation}
$\partial^C\varphi (\bar{u})$ reduces to a singleton and equality holds in \eqref{clarkeinclu}, if additionally, the condition
\begin{equation}\label{measurezero2}
\mu_\zeta (\{v\in\mathbb{S}^{m-1}\mid\# \partial^C_ue(\bar{u},v)\geq 2\})=0
\end{equation}
is satisfied. As a consequence, $\varphi$ is strictly differentiable at $\bar{u}$ in the Hadamard sense \cite[p.~30]{clarke1983}.
\end{theorem}

\noindent We note that differentiability of $\varphi$, and thus condition \eqref{measurezero2}, may be violated in general (see \cite[Example 2.15]{geiersbach2023optimality}). Therefore, the question arises if there are weaker conditions than those enforcing differentiability as in \eqref{measurezero2}, which still guarantee equality in \eqref{clarkeinclu}. It was shown in \cite[Theorem 2.3]{geiersbach2023optimality} that this holds true under the additional assumptions that $g$ is jointly (!) convex in both variables and that the set 
\begin{equation}\label{boundedset}
\{z \in \R^m \mid g(\bar{u},z)\leq 0\}
\end{equation}
is bounded. The condition \eqref{boundedset} turns out to be quite restrictive (for instances where it is violated, see \cite[Examples 2.15 and 2.16]{geiersbach2023optimality}) and it is difficult to check, in general. Now, we are going to state a stronger result, which allows us to avoid the aforementioned conditions at the price of a growth condition which, however, turns out to be automatically satisfied in our control problem \eqref{eq:probuniform-problem}. We shift the rather lengthy proof to the appendix.
\begin{theorem}\label{clarkeequality}
Let some $\bar{u}\in U$ be given. For $(u,z)\in U\times\mathbb{R}^m$ assume that $g(u,z)=\sup_{w\in K}\,h(u,z,w)$, where $K\subseteq\mathbb{R}^l$ is compact and $h\colon U\times\mathbb{R}^m\times K\to\mathbb{R}$ satisfies the conditions
\begin{enumerate} 
\item $h$ is continuous, Fréchet differentiable in its first two arguments, and convex in its second argument;
\item the mapping $(u,z,w) \mapsto (D_{u} h(u,z,w),D_z h(u,z,w) )$ is continuous; 
\item $h(\bar{u},0,w)<0\quad\forall w\in K$;
\item $\exists c>0:\left\Vert D_{u}h(u,z,w) \right\Vert \leq c\exp\left({\left\Vert
z\right\Vert }\right)\,\,\forall u:\left\Vert u-\bar{u}\right\Vert \leq
c^{-1}\,\,\forall z:\|z\|\geq c\,\,\forall w\in K $. 
\end{enumerate}
Then, $-\varphi$ is regular at $\bar{u}$ in the sense of Clarke and equality holds in \eqref{clarkeinclu}.
\end{theorem}

\subsection{(Sub-)differential of the probability function in the concrete optimal control problem}
We now want to calculate the exact subdifferential of the probability function associated with problem \eqref{eq:probuniform-problem}. A corresponding result has been previously obtained in \cite[Theorem 2.13]{geiersbach2023optimality}. However, the restrictive condition requiring the set \eqref{boundedset} to be bounded was imposed there. The results of the previous section allow us to get rid of this assumption, which we will now show. To this end, denote by $S\colon L^2(D) \times \mathbb{R}^m \rightarrow \mathcal{C}(\bar{D})$ the parameterized control-to-state operator assigning to each $u\in L^2(D)$ and each $z\in\mathbb{R}^m$ the solution $y$ of the PDE \eqref{eq:probuniform-problem-b}-\eqref{eq:probuniform-problem-c} with right-hand side $u+f(\cdot,z)$. Due to Assumption~\ref{ass:PDE-standing} and \cite[Lemma 1.2]{geiersbach2023optimality}, the operator $S$ is well-defined and bounded.  Then the control problem \eqref{eq:probuniform-problem} may be recast in the general form of \eqref{optprob} upon defining the function $g(u,z):=\sup_{x\in D}[S(u,z)](x)-\alpha$. Since $S$ maps into $\mathcal{C}(\bar{D})$, we may write
\begin{equation}\label{gdef}
g(u,z):=\max\limits_{x\in\bar{D}}[S(u,z)](x)-\alpha .  
\end{equation}
Due to linearity, one has the decomposition
\begin{equation}\label{sstructure}
S(u,z)=\bar{y}_{(u)}+\sum\limits_{i=1}^mz_iy^{(i)}\quad (u\in L^2(D), z\in\mathbb{R}^m)
\end{equation}
where $\bar{y}_{(u)}$ is the solution of the PDE
\begin{equation}\label{meanstate}
-\Delta  y(x)=u(x)+f_0(x)\quad x \in D; \quad y(x)=0,\quad x \in \partial D
\end{equation}
and the $y^{(i)}$ are the solutions of the PDEs 
\begin{equation}\label{basicstates}
-\Delta  y(x)=\phi_i(x)\quad x \in D; \quad y(x)=0,\quad x \in \partial D  \quad (i=1,\ldots, m).
\end{equation}
We shall refer to $\bar{y}_{(u)}$ as the mean state associated with the control $u$ (because it relates to the mean value $f_0$ of the random source term) and to the $y^{(i)}$ as the basic states (because they relate to the basic functions $\phi_i$ in the random source term with no control in action).
Moreover, to each $x\in D$ we assign a dual element ${\bf u}_x\in (L^2(D))^*$ by ${\bf u}_x(h):={y^{h}(x)}$ for all $h\in L^2(D)$, where $y^{h}$ is the solution of the control-only PDE
\begin{equation}\label{controlonly}
-\Delta  y(x)=h(x)\quad x \in D; \quad y(x)=0,\quad x \in \partial D.
\end{equation}
\begin{theorem}\label{exactclarke}
In the control problem \eqref{eq:probuniform-problem}, fix some point of interest $\bar{u}\in L^2(D)$. Assume that the mean state $\bar{y}_{(\bar{u})}$ associated with $\bar{u}$ satisfies the condition
\begin{equation}\label{sp1}
\bar{y}_{(\bar{u})}(x)<\alpha\quad\forall x\in\bar{D}.
\end{equation}
Then, the probability function $\varphi$ associated with our control problem is locally Lipschitzian at $\bar{u}$ and has the exact Clarke subdifferential 
\[
\partial^C\varphi (\bar{u})=-\int\limits_{\{v\in\mathbb{S}^{m-1}:\rho (v)<\infty\}}{\rm clco}\left\{\frac{f_\chi (\rho (v))}{\varkappa (v,x)}\cdot {\bf u}_x\mid x\in M^*(v)\right\} \D\mu_\zeta (v),
\]
where ``clco'' means the closed convex hull and for each $v\in\mathbb{S}^{m-1}$ and $x\in\bar{D}$, 
\begin{eqnarray*}
\varkappa (v,x)&:=&\sum\limits_{i=1}^m(\Sigma^{1/2}v)_iy^{(i)}(x),\\ \rho (v)&:=&\max\{r\geq 0\mid 
\bar{y}_{(\bar{u})}(x)+r\varkappa (v,x)\leq\alpha\,\,\forall x\in\bar{D}\},\\
M^*(v)&:=&\{x\in\bar{D}\mid\varkappa (v,x)>0,\,\,
\bar{y}_{(\bar{u})}(x)+\rho (v)\varkappa (v,x)=\alpha\}.
\end{eqnarray*}
\end{theorem}
\proof It has been shown in \cite[Lemma 2.11, Lemma 2.12]{geiersbach2023optimality} that the integral 
\[
\int\limits_{\mathbb{S}^{m-1}}\partial^C_ue(\bar{u},v) \D \mu_\zeta (v)
\]
coincides with the right-hand side of the identity in the statement of this theorem. Therefore, it will be sufficient to apply Theorem \ref{clarkeequality} to our control problem \eqref{eq:probuniform-problem} in order to use \eqref{clarkeinclu} as an equality. Clearly, this would yield the asserted exact formula. Hence, we only have to check the assumptions of Theorem \ref{clarkeequality}. Given the definition \eqref{gdef}, we set $K:=\bar{D}$, $w:=x$, $h(u,z,x):=[S(u,z)](x)-\alpha$ and observe that we are in the setting of Theorem \ref{clarkeequality}. We check the assumptions of that theorem: as for assumption 1., let us first verify the continuity of $h$ by considering a sequence $(u_n,z_n,x_n)\to (u,z,x)$ in the set $L^2(D)\times\mathbb{R}^m\times\bar{D}$. 
As shown as part of the proof in \cite[Lemma 2.7]{geiersbach2023optimality}, there exists a constant $C>0$ such that for all $(u_1,z_1),(u_2,z_2)\in L^2(D)\times\mathbb{R}^m$,
\begin{equation}\label{lipschitz}
\|[S(u_1,z_1)]-[S(u_2,z_2)]\|_{\mathcal{C}(\bar{D})}\leq C (\|u_1-u_2\|_{L^2(D)}+\|z_1-z_2\|).
\end{equation}
Accordingly, 
\begin{eqnarray*}
&|h(u_n,z_n,x_n)-h(u,z,x)|\leq&\\&|[S(u_n,z_n)](x_n)-[S(u,z)](x_n)|+|[S(u,z)](x_n)-[S(u,z)](x)|\leq&\\&C(\|u_n-u\|_{L^2(D)}+\|z_n-z\|)+|[S(u,z)](x_n)-[S(u,z)](x)|&
\end{eqnarray*}
Here, both terms on the right converge to zero (the second one thanks to $S(u,z)\in\mathcal{C}(\bar{D})$). This shows the continuity of $h$. As observed in \cite[p.3]{geiersbach2023optimality}, the operator $S$ admits for some $y_0\in\mathcal{C}(\bar{D})$ a decomposition
\[
S(u,z)=P(u,z)+y_0\quad\forall (u,z)\in L^2(D) \times \mathbb{R}^m 
\]
where $P\colon L^2(D) \times \mathbb{R}^m\to\mathcal{C}(\bar{D})$ is some continuous linear operator. Consequently, for fixed $x\in\bar{D}$, the function $(u,z)\mapsto h(u,z,x)$ is the sum of a constant $y_0(x)-\alpha$ and a continuous linear function $(u,z)\mapsto [P(u,z)](x)$. This implies that $h$ is Fréchet differentiable and convex in its first two arguments. Moreover, the partial Fréchet derivative $D_{(u,z)}h$ at $(u,z,x)$ equals the continuous linear function $[P(\cdot,\cdot)](x)$. Altogether, we have shown the validity of assumptions 1. and 2. of Theorem \ref{clarkeequality}. Assumption 3. follows immediately from \eqref{sp1}. In order to check assumption 4., observe that by \eqref{lipschitz}, the function $u\mapsto h(u,z,x)$ is locally Lipschitzian with some common modulus $C$ at all $u\in L^2(D)$ and for all $(z,x)\in\mathbb{R}^m\times\bar{D}$. Consequently,
\[
\left\Vert D_{u}h(u,z,x) \right\Vert\leq C\leq ce^{\|z\|}\,\,\forall (u,z,x)\in L^2(D)\times\mathbb{R}^m\times\bar{D}:\|z\|\geq c:=\max\{1,\log C\}.
\]
This shows the validity of assumption 4. and finishes the proof.

\qed
\begin{remark}
One may easily derive from \cite[Prop. 3.11]{vanAckooij_Henrion_2014} that condition \eqref{sp1} is satisfied whenever the probability $\varphi (\bar{u})$ is larger than or equal to 0.5. Since typically the probability level $p$ is chosen close to one, this condition will be automatically fulfilled for most iterations of the control. Hence, the exact formula for the Clarke subdifferential in Theorem \ref{clarkeequality} will basically hold true unconditionally. 
\end{remark}
\begin{remark}
If the sets $M^*(v)$ introduced in Theorem \ref{exactclarke}   satisfy the condition
\begin{equation}\label{singleton}
\#M^*(v)=1\quad\mu_\zeta-{\rm a.e.}\quad v\in\mathbb{S}^{m-1},
\end{equation}
then the integral in the formula for the Clarke subdifferential, and hence the Clarke subdifferential itself, reduce to singletons. Therefore, the probability function is strictly differentiable in the Hadamard sense (see \cite[Prop. 2.2.4]{clarke1983}) and its derivative is given by 
\begin{equation}\label{strictderiv}
D\varphi (\bar{u})=-\int\limits_{\{v\in\mathbb{S}^{m-1}:\rho (v)<\infty\}}\frac{f_\chi (\rho (v))}{\varkappa (v,x^*(v))}\cdot {\bf u}_{x^*(v)} \D \mu_\zeta (v),
\end{equation}
where, for $\mu_\zeta-$ almost every $v\in\mathbb{S}^{m-1}$,  $x^*(v)$ is defined as the unique element in the set $M^*(v)$.
\end{remark}
Simple examples show that in the setting of Theorem \ref{exactclarke} the probability function may fail to be differentiable, hence condition \eqref{singleton} is violated in general.
Evidently, \eqref{singleton} is difficult to verify for concrete data. 
In \cite[Lemma 4.3]{ackooij-henrion2017}, some easy to understand constraint qualification (so-called ``rank 2-CQ'') has been shown to imply \eqref{singleton} for finite random inequality systems. A generalization to our setting with infinite systems (uniform state constraints) does not seem to be straightforward, so that the verification of differentiability for the probability function remains an open problem. On the other hand, in the numerical solution of the given control problem, one may be forced to replace the uniform state constraint \eqref{eq:probuniform-problem-d} by evaluating it on a finite subset $\tilde{D}\subseteq D$ of the domain. Then, the aforementioned rank 2-CQ reduces to the following verifiable condition at some fixed control $\bar{u}\in U$:
\begin{eqnarray*}
&{\rm rank}\,\{Y(x^a),Y(x^b)\}=2\quad\forall x^a,x^b\in\tilde{D}:&\\ x^a\neq x^b,\quad 
&\bar{y}_{(\bar{u})}(x^a)+\langle z,Y(x^a)\rangle =\bar{y}_{(\bar{u})}(x^ b)+\langle z,Y(x^b)\rangle =\alpha&\\&\forall z\in\mathbb{R}^m: \bar{y}_{(\bar{u})}(x)+\langle z,Y(x)\leq\alpha\,\,\forall x\in\tilde{D}&
\end{eqnarray*}
Here, for $x\in\tilde{D}$, the vector $Y(x):=(y^{(1)}(x),\ldots,y^{(m)}(x))$ collects in its components all values of the basic functions $y^{(i)}$ at $x$.
\subsection{Convexity properties under Gaussian and truncated Gaussian distributions}\label{convtrunc}

\noindent
In this section, we show that problem \eqref{eq:probuniform-problem} is convex under Gaussian and truncated (to convex sets) Gaussian distributions. Given the assumed convexity of the objective $F$, it is sufficient to verify convexity of the constraint. For pure Gaussian distributions, this has been done (without explicit statement) in \cite{geiersbach2023optimality}. We refer to the abstract form \eqref{optprob} of problem \eqref{eq:probuniform-problem}, where $\varphi$ is defined as the Gaussian probability function associated with the random inequality $g(u,\xi(\omega))\leq 0$. In the context of our problem \eqref{eq:probuniform-problem}, this constraint function takes the form \eqref{gdef}. Similarly, we may consider problem \eqref{eq:probuniform-problem} under a truncated Gaussian distribution and represent it in the abstract form \eqref{optprob} with $\varphi$ being replaced by the truncated Gaussian probability function $\tilde{\varphi}$
associated with the same constraint function $g$ from \eqref{gdef}. More precisely, if $\tilde{\xi}\sim\mathcal{TN}(\mu,\Sigma,\Xi)$ designates a random vector having a Gaussian distribution $\mathcal{N}(\mu,\Sigma)$ truncated to a closed set $\Xi:=\{z\in\mathbb{R}^m\mid s(z)\leq 0\}$ ($s\colon\mathbb{R}^m\to \mathbb{R}$ continuous), then
\[
\tilde{\varphi}(u)=\mathbb{P}(\omega\mid g(u,\tilde{\xi}(\omega))\leq 0)=\frac{\mathbb{P}(\omega\mid g(u,\xi (\omega))\leq 0,\,\,s(\xi (\omega))\leq 0)}{\mathbb{P}(\omega\mid s(\xi (\omega))\leq 0)}.
\]
The last equation above expresses the fact that the truncated distribution is nothing but the original distribution conditioned to a subset.
Consequently, the probabilistic constraint $\tilde{\varphi} (u)\geq p$ in \eqref{optprob} (with $\varphi$ replaced by $\tilde{\varphi}$) can be equivalently written as $\mathbb{P}(\omega\mid\tilde{g}(u,\xi(\omega))\leq 0)\geq\tilde{p}$, where 
\begin{equation}\label{gtildedef}
\tilde{g}(u,z):=\max\{g(u,z),s(z)\},\quad\tilde{p}:=p\cdot\mathbb{P}(\omega\mid s(\xi (\omega))\leq 0).
\end{equation}
\begin{proposition}\label{convexity}
The optimization problem \eqref{eq:probuniform-problem} is convex if the random vector $\xi$ follows a Gaussian distribution or a Gaussian distribution truncated to a set
$\Xi:=\{z\in\mathbb{R}^m\mid s(z)\leq 0\}$  where $s\colon \mathbb{R}^m\to \mathbb{R}$ is convex, hence continuous.
\end{proposition}
\proof
The convexity of the constraint $\varphi(u)\geq p$ in case of a Gaussian distribution follows in the abstract setting from the convexity of the constraint function $g(u,z)$ (jointly in both variables), see \cite[Lemma 2.4]{geiersbach2023optimality}. That this property holds true for the concrete function $g$ defined in \eqref{gdef}, has been shown in \cite[Lemma 2.7]{geiersbach2023optimality}. As for the case of a truncated distribution, we have seen above that we may represent the corresponding chance constraint as in the Gaussian case but with the modified constraint function $\tilde{g}$ and the modified probabilty level $\tilde{p}$ defined in \eqref{gtildedef}. Evidently, $\tilde{g}$ is jointly convex in $(u,z)$ since $g$ is so and $s$ is convex. Consequently, we can apply again the previous argument to prove the convexity of the feasible set also under truncated Gaussian distribution.
\qed
\section{Numerical results for the control problem with probabilistic constraints}\label{numprobcons}
\subsection{Implementation of the spherical-radial decomposition}\label{spherapprox}
In order to deal with the probabilistic constraint $\varphi(u)\geq p$ in the general problem \eqref{optprob} numerically, one has to appropriately approximate the probability function $\varphi$ and, assuming differentiability, its derivative  $D\varphi$. Since $\varphi$ is defined as the spherical integral \eqref{srd}, it can be approximated by a finite sum 
\[
\varphi (u)\approx\hat{\varphi}(u):=\frac{1}{K}\sum_{k=1}^Ke(u,v^{(k)})\quad\forall u\in U,
\]
where $\{v^{(k)}\}_{k=1}^K$ is a sample of size $K$ of the uniform distribution on the sphere $\mathbb{S}^{m-1}$. A simple way to create such a sample consists in normalizing to unit length a quasi-Monte Carlo (QMC) sample of the given Gaussian distribution. It was shown in \cite[eq. (19)]{geiersbach2023optimality} that, for any fixed $\bar{u}$ with $g(\bar{u},0)<0$, there exists a neighborhood $\mathcal{N}$ of $\bar{u}$ such that 
\begin{equation}\label{erep}
e(\bar{u},v)=\left\{
\begin{array}{ll}
F_\chi(\rho^*(v))&\mbox{if }\rho^*(v)<\infty\\
1&\mbox{if }\rho^* (v)=\infty
\end{array}\right.\quad\forall u\in\mathcal{N}\,\,\forall v\in\mathbb{S}^{m-1},
\end{equation}
where $F_\chi$ is the cumulative Chi-distribution function introduced before  and 
\[
\rho^* (v):=\sup\{r\geq 0\mid g(\bar{u},r\Sigma^{1/2}v)\leq 0\}\quad (v\in\mathbb{S}^{m-1}).
\]
It follows from \cite[eq. (29)]{geiersbach2023optimality}  that $\rho^*$ coincides with the function $\rho$ defined in the statement of Theorem \ref{exactclarke}. Consequently, the probability function can be approximated at some $\bar{u}$ with $g(\bar{u},0)<0$ by
\begin{eqnarray*}
\varphi (\bar{u})\approx\hat{\varphi}(u)&=&\frac{1}{K}\#\{k\in\{1,\ldots ,K\}\mid \rho (v^{(k)})=\infty\}\\&&+\frac{1}{K}\sum_{\{k\in\{1,\ldots ,K\}:\rho (v^{(k)})<\infty\}}F_\chi (\rho (v^{(k)})).
\end{eqnarray*}
Should $\varphi$ be differentiable, the spherical integral in \eqref{strictderiv} can be approximated by the sample-based quantity
\[
D\varphi (\bar{u})\approx D\hat{\varphi} (\bar{u}):=-\sum_{\{k\in\{1,\ldots ,K\}:\rho (v^{(k)})<\infty\}}\frac{f_\chi (\rho (v^{(k)}))}{\varkappa (v^{(k)},x^*(v^{(k)}))}\cdot {\bf u}_{x^*(v^{(k)})},
\]
whose ingredients are defined in the statement of Theorem \ref{exactclarke} and below \eqref{strictderiv}, respectively. 

It is essential to observe that values and derivatives (in a specified direction $h\in L^2(D)$) can be simultaneously updated by sharing the same sample $v^{(k)}$. In this way, the potentially time-consuming computation of $\rho (v^{(k)})$ may be carried out only once, when it comes to updating the probabilities themselves, whereas it does not have to be redone during the update of the gradient.
Altogether, this leads us to the algorithmic scheme in Algorithm \ref{alg0}.
\begin{algorithm}[htb] 
\caption{Algorithmic approximation of probabilities and their gradients via spherical-radial decomposition}\label{alg0} 
{\bf Input}: $D\subset\mathbb{R}^d$ (domain), $f_0$ (expected source term), $\phi_i (i=1,\ldots m)$ (basic source functions), $\Sigma$ (covariance matrix of Gaussian distribution for random coefficients), $\bar{u}$ (control of interest satisfying \eqref{sp1}), $h$ (direction of interest to apply the gradient to).
\begin{enumerate}
\item Initialization
\begin{enumerate}
\item Fix sample size $K$ and create a sample 
$\{v^{(k)}\}_{k=1}^K$ uniformly distributed on the sphere $\mathbb{S}^{m-1}$.
\item $prob:=0, deriv:=0, k:=1$.
\item Determine functions $y^{(i)} (i=1,\ldots m)$ as the solutions of \eqref{basicstates}.
\item Determine function $\bar{y}_{(\bar{u})}$ as the solution to \eqref{meanstate} (with $u:=\bar{u}$).
\item Determine function $y^h$ as the solution to \eqref{controlonly}.
\end{enumerate}
\item Iteration on $k$ 
\begin{enumerate}
\item Define a function $\varkappa:=\sum\limits_{i=1}^m(\Sigma^{1/2}v^{(k)})_iy^{(i)}$
\item Calculate $\rho:=\max\{r\geq 0\mid 
\bar{y}_{(\bar{u})}(x)+r\varkappa (x)\leq\alpha\,\,\forall x\in\bar{D}\}$
\item Determine a point $x^*\in\bar{D}$ satisfying
$\bar{y}_{(\bar{u})}(x^*)+\rho\varkappa (x^*)=\alpha$ 

(if $x^*$ is not unique then the probability function fails to be differentiable and the output will yield a subgradient rather than a gradient).
\item If $\rho<\infty$ then $prob:=prob+F_\chi (\rho)$ else $prob:=prob+1$.
\item If $\rho<\infty$ then $deriv:=deriv+\frac{f_\chi (\rho)}{\varkappa(x^*)}y^h(x^*)$.
\end{enumerate}
\item Termination: If $k<K$ then $k:=k+1$. Goto Step 2.
\item Output: $prob:=K^{-1}prob, deriv:=K^{-1}deriv$
\end{enumerate}
{\bf Output:} Approximation of probability $\varphi (\bar{u})\approx prob$ and of its directional derivative $D\varphi (\bar{u})(h)\rangle\approx deriv$.
\end{algorithm}
\subsection{Results in dimension one}\label{Sect_resulstindimone}
We consider problem \eqref{eq:probuniform-problem} with the following data:
\begin{small}
\begin{align*}
&d=1, D=(0,1), m=6, \alpha=0.2, F(u)=\|u\|^2_{L^2(D)}, p=0.9, \Sigma_{i,j}=9\cdot 0.6^{|i-j|}\,\, (i,j=1,\ldots m),\\ &f_0(x)=5x^2, 
\phi_{2i-1}(x)=\sin (ix)\,\, (i=1,2,3), 
\phi_{2i-2}(x)=\cos (x/i)\,\, (i=2,3,4).
\end{align*}
\end{small}

\begin{figure}[htb]
    \centering
\includegraphics[width=7cm]{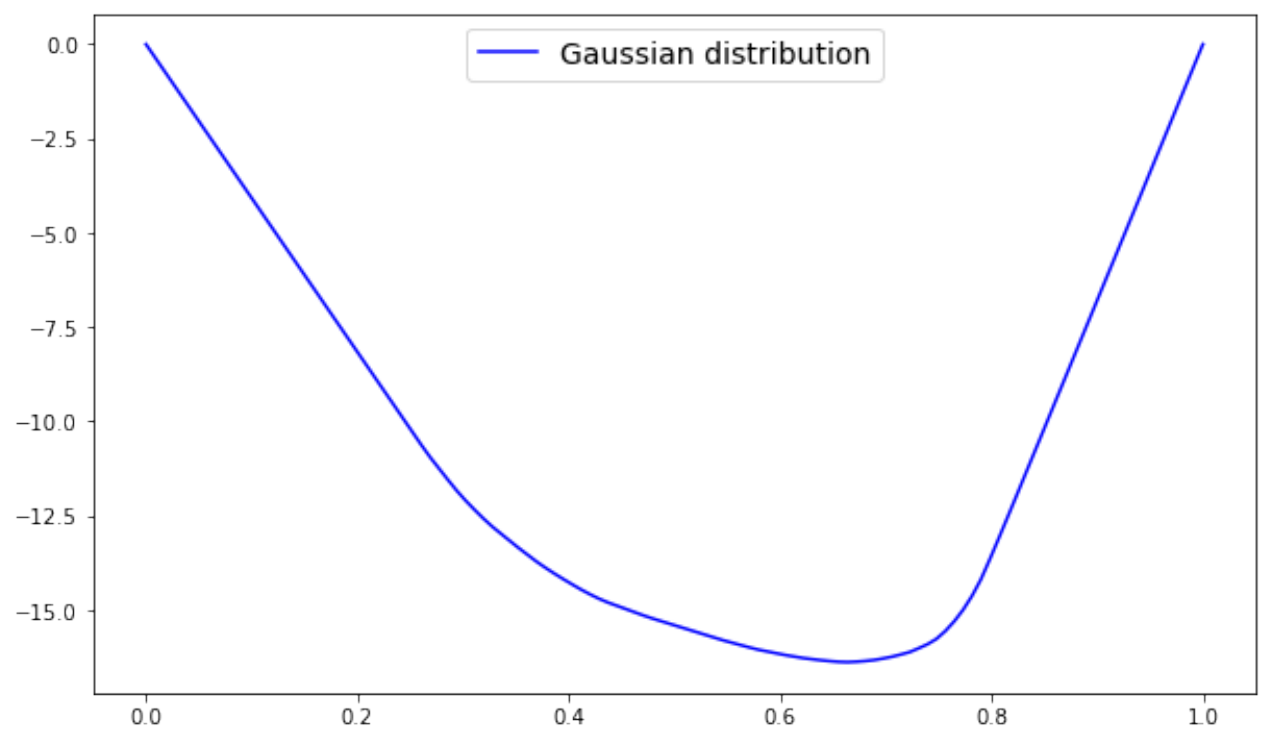}
\includegraphics[width=7cm]{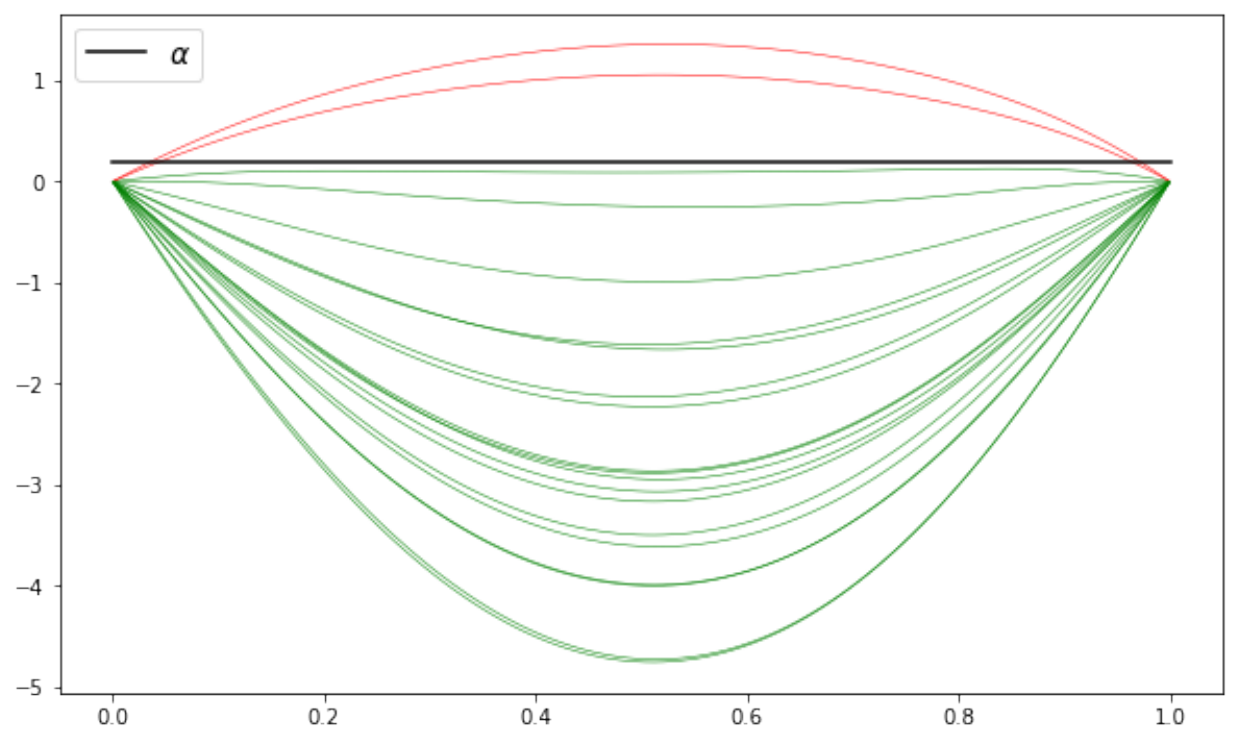}
    \caption{Numerical solution of \eqref{eq:probuniform-problem}. Optimal control for probability level $p=0.9$ (left) and associated states for twenty scenarios of the random source term (right)}
    \label{fig1}
\end{figure}

\noindent
We use a finite difference discretization over a subdivision of the domain into 120 intervals. The values and---assuming differentiability---derivatives of the probability function in were obtained as described in Algorithm \ref{alg0} on the basis of a QMC sample on the unit sphere of size 512. The optimization problem was numerically solved with the ``SLSQP'' method from the Python standard minimzation routine {\it scipy.optimize.minimize}. Fig. \ref{fig1} (left) shows the optimal control under the indicated Gaussian distribution.
Fig. \ref{fig1} (middle) plots twenty states associated with the optimal control and with twenty scenarios for the random source term $f(x,\xi (\omega))$ in \eqref{eq:probuniform-problem-b}. It can be seen that only two of them exceed the desired threshold $\alpha$ occasionally, whereas the remaining ones stay below it uniformly over the domain, which corresponds well to the imposed probability level of $p=0.9$. 
\begin{figure}[htb]
    \centering
\includegraphics[width=7cm]{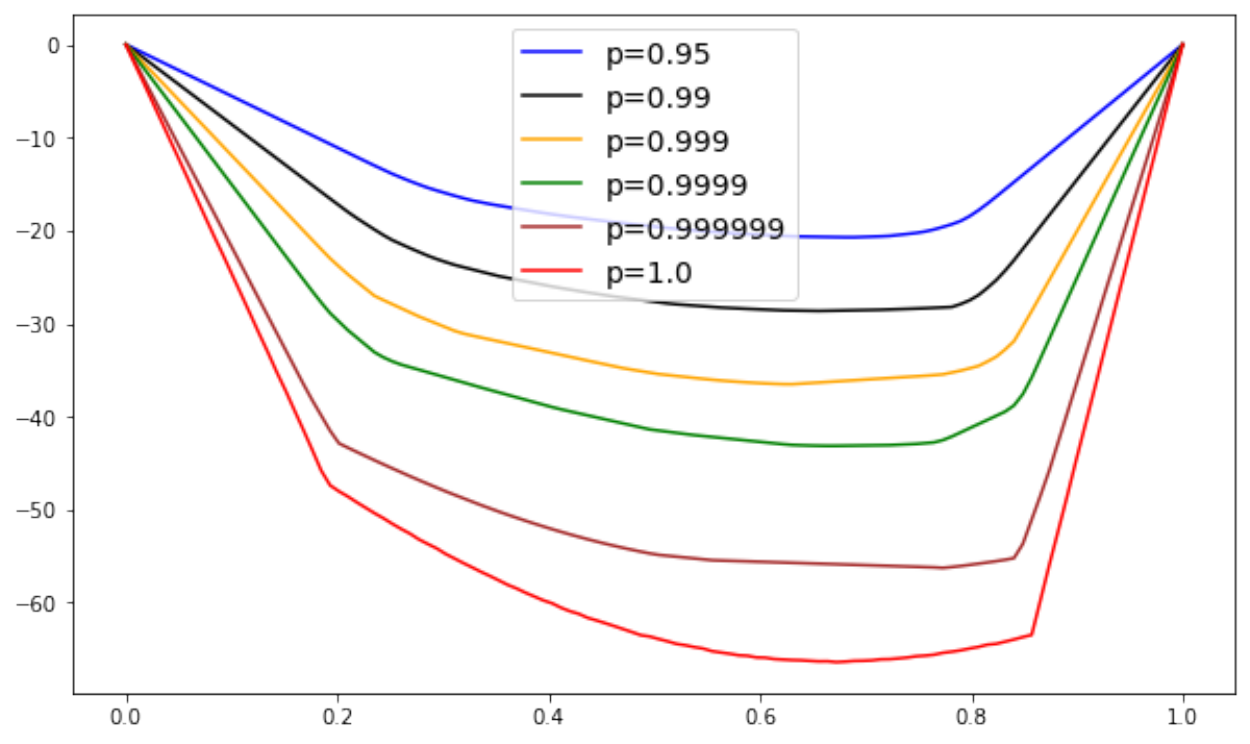}
    \caption{Solution of \eqref{eq:probuniform-problem} for Gaussian distribution truncated to an ellipsoid and for increasing probability levels.}
    \label{fig1a}
\end{figure}

In order to include the almost sure case ($p=1$) for later comparison, the problem was solved again with the Gaussian distribution truncated to an ellipsoid $\mathcal{E}$ defined by $\mathcal{E}:=\{z\in\mathbb{R}^m\mid z^T\Sigma^{-1}z\leq 36\}$ (the choice of the Mahalanobis norm  in this definition facilitates computations). In this way, the support of the distribution becomes compact so that the problem with an almost-sure constraint has a chance to have a feasible solution. As can be seen from Fig. \ref{fig1} (right), the optimal solutions of the probabilistically constrained problems seem to converge to that with almost sure constraints for $p\to 1$ (for a rigorous statement, see Section \ref{aslimitpc}). 

Note that passing to truncated Gaussian distributions already points to the possibility of using alternative distributions. We do not present details here but note that the same methodology could be applied to multivariate lognormal or Student or Gaussian mixture distributions.

\subsection{Results in dimension two}
In this section, we revisit problem \eqref{eq:probuniform-problem} but with some increased complexity when compared to the previous section: the domain will be in two dimensions, the dimension of the random parameter will increase from 6 to 30, and additional bound constraints will be imposed on the control. More precisely, we consider
the following data:

\begin{small}
\begin{align*}
&d=2, D=(0,1)^2, m=30, \alpha=0.2, F(u)=\|u\|^2_{L^2(D)}, p=0.9, \Sigma_{i,j}=9\cdot 0.6^{|i-j|}\,\, (i,j=1,\ldots m),\\ &f_0(x_1)=5x_1x_2,
 \phi_i(x_1,x_2)=\sin (ix_1)\cos (ix_2)\,\, (i=1,\ldots ,m).
\end{align*}
\end{small}
In addition, we impose the following bounds on the control: 
\[
-5\leq u(x_1,x_2)\leq 0\quad {\rm a.e.} \,(x_1,x_2)\in D.
\]
\begin{figure}[htb]
    \centering
    \includegraphics[width=7cm]{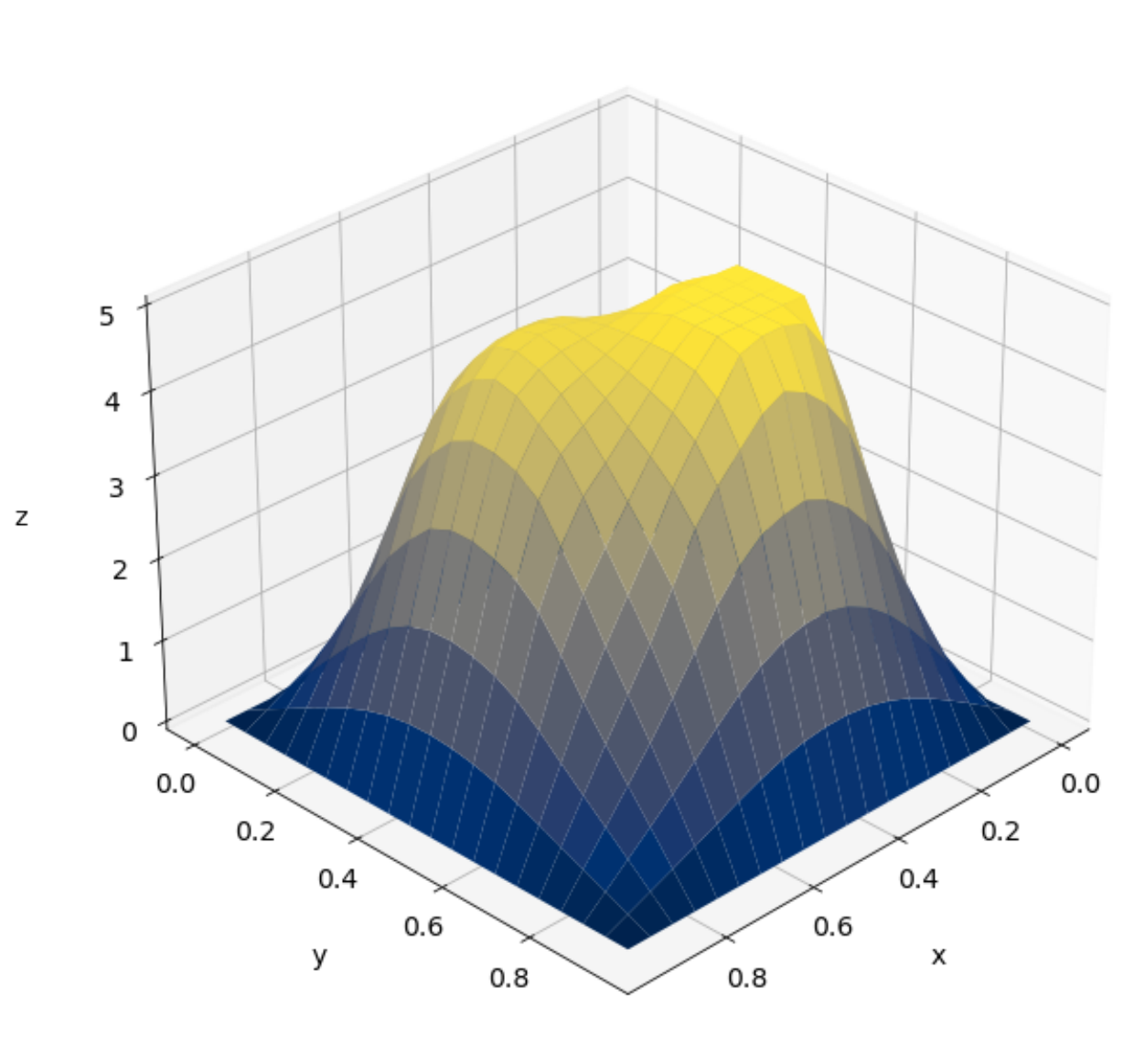}
    \hspace{0.5cm}
    \includegraphics[width=6.6cm,height=6cm]{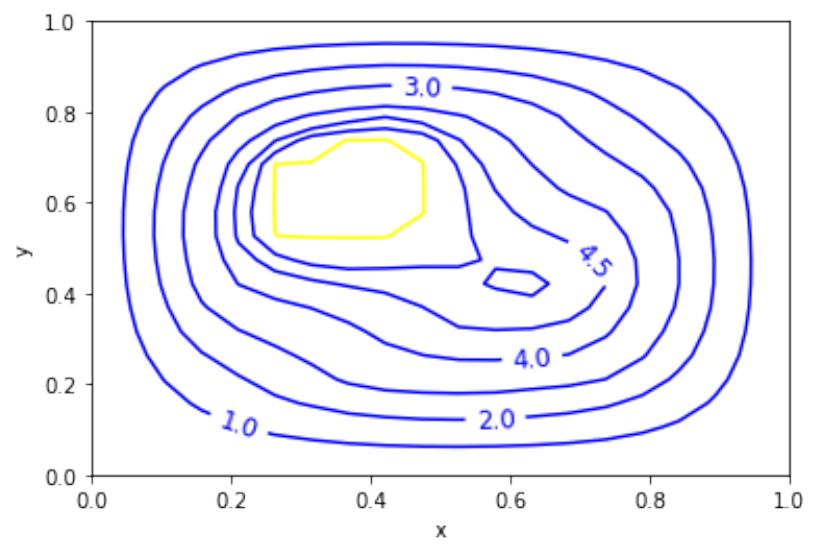}
    \caption{Optimal solution of the control problem in 2D (left). For graphical reasons the negative control is plotted. Contourplot (right).}
    \label{solution_2d}
\end{figure}
The PDE is numerically solved using a Python code by Zaman available on 
\begin{center}
{\tt https://github.com/zaman13/Poisson-solver-2D}   
\end{center}
and which is based on the paper \cite{zaman2022}. The PDE was solved using finite differences on a 20x20 grid of the domain and the random parameter was approximated by using $2^{13}=8192$ QMC samples on the sphere.
Figure \ref{solution_2d} shows the resulting optimal control (linearly interpolated). It can be seen in both plots that the control bound becomes active in a certain small region of the domain (yellow).
\section{Almost sure constraints}
\label{sec:almost-sure}
In this section, we want to consider the extreme case of risk-averse decision making, namely the optimal control with constraints that should be satisfied almost surely. In this case, our PDE-constrained optimization problem becomes 
\begin{subequations}
    \label{eq:a.s.-problem}
    \begin{alignat}{3}
    \min_{u \in L^2(D)}\, F(u) & & &&  \\
   \text{s.t.}  \quad -\Delta  y(x,\omega)  &=  u(x) + f(x,\xi (\omega)), & &&\quad x \in D \quad \text{$\mathbb{P}$-a.s.}, \label{eq:a.s.-problem-constraint1} \\
 y(x,\omega) &=0, & &&\quad x \in \partial D  \quad \text{$\mathbb{P}$-a.s.}, \label{eq:a.s.-problem-constraint2}\\
 y(x,\omega) &\leq \alpha, & && x\in D \quad \text{$\mathbb{P}$-a.s.}\phantom{.} \label{eq:a.s.-problem-constraint3}
    \end{alignat}
\end{subequations}
We shall assume in this section that the support $\Xi\subseteq\mathbb{R}^m$ of the random vector $\xi$ is compact, because otherwise there may not exist a feasible solution that satisfies \eqref{eq:a.s.-problem-constraint3}. This assumption will meet the setting illustrated in Fig.~\ref{fig1} (right), where the distribution was assumed to be Gaussian, truncated to an ellipsoid. 

\subsection{Almost sure constraints as limit of probabilistic constraints}\label{aslimitpc}
Evidently, by equivalence of \eqref{eq:a.s.-problem-constraint3} with \eqref{eq:probuniform-problem-d} when choosing $p=1$, this optimization is equivalent with the previously analyzed probabilistically constrained problem \eqref{eq:probuniform-problem} when choosing the maximum possible probability level. We emphasize that, in spite of this equivalence, necessary and sufficient optimality conditions can be obtained for 
the formulation \eqref{eq:a.s.-problem} (\cite[Theorem 3.10]{geiersbach2023optimality}) while they cannot for the formulation \eqref{eq:probuniform-problem} (\cite[Remark 2.6 and Example 3.1]{geiersbach2023optimality}). On the other hand, we were able numerically to find a candidate solution to \eqref{eq:probuniform-problem} for $p=1$ and observed a convergence of solutions for $p\to 1$ towards this candidate (see Fig. \ref{fig1a}).
Two questions arise: is the candidate we found the true solution of the almost sure problem \eqref{eq:a.s.-problem} and can we prove the mentioned convergence? The first question will be considered in the next section, whereas the second question will be answered after the following technical preparation.
\begin{lemma}\label{prep}
Consider problem \eqref{eq:probuniform-problem} but with the random vector $\xi$ having a Gaussian distribution truncated to a compact convex set $\Xi$ as in Proposition \ref{convexity}.
Let   $u_{p_k}\in L^2(D)$ be optimal solutions of \eqref{eq:probuniform-problem}  with $ p_k \in (0,1)  $ and $p_k \to 1$. In addition, suppose that $ u_{p_k} \rightharpoonup  \bar{u} \in L^2(D)$. Then, $\bar{u}$ is an optimal solution of problem \eqref{eq:probuniform-problem} with $p=1$ and $F(u_{p_k}) \to F(\bar{u})$.    
\end{lemma}
\proof
As shown in Section \ref{convtrunc}, the constraint of 
\eqref{eq:probuniform-problem} can be rewritten as $\tilde{\varphi} (u)\geq\tilde{p}$, where $\tilde{\varphi} (u)=\mathbb{P}(\omega\mid\tilde{g}(u,\xi(\omega))\leq 0)$ is a probability function related to an untruncated Gaussian random vector and $\tilde{g}, \tilde{p}$ are defined in \eqref{gtildedef}.
By \cite[eq. (5)]{geiersbach2023optimality}, the functions $(u,z)\mapsto [S(u,z)](x)$  are affine linear (hence weakly lower semicontinuous) for each $x\in\bar{D}$. Consequently, the function $g$ in \eqref{gdef} is weakly lower semicontinuous as a sum of a constant and a maximum of such functions. On the other hand, the function $(u,z)\mapsto s(z)$  with $s$ from Proposition \ref{convexity} is convex and continuous, hence weakly lower semicontinuous. It follows once more that $\tilde{g}$ as the maximum of two weakly lower semicontinuous functions shares this property.  
Therefore, using \cite[Lemma 2]{Farshbaf-Shaker2018},
we get that the function $\varphi$ defined in \eqref{optprob} is weakly sequentially upper semicontinuous, whence
$\varphi(\bar{u}) \geq \limsup_{k \to \infty  } \varphi(u_{p_k}) = 1$.
Now, let $u' \in L^2(D)$ be any   control function satisfying $\varphi(u') = 1$. Since in problem \eqref{optprob} with probability level $p_k$, the control $u'$ is feasible while the control $u_{p_k}$ is optimal, it follows that $F(u_{p_k}) \leq F(u')$. In particular, it shows that $F(u_{p_k}) \leq F(\bar{u})$ and consequently $\limsup_{k \to \infty}F(u_{p_k}) \leq F(\bar{u})$. Therefore, utilizing the weak lower semicontinuity of the objective function, we can conclude that $F(\bar{u}) \leq \liminf_{k\to \infty} F(u_{p_k}) \leq F(u')$, thereby demonstrating the optimality of $\bar{u}$ and $ \lim_{k\to \infty} F(u_{p_k}) = F(\bar{u})$. 
\qed

\begin{proposition} \label{propoPartial}
In addition to the assumptions of Lemma \ref{prep}, assume that the objective $F$ in \eqref{eq:probuniform-problem} is strictly convex. Then, any bounded sequence $(u_{p_k})_{k\in  \mathbb{N}}$ of minimizers of \eqref{eq:probuniform-problem}  with $p_k \in (0,1)$  weakly converges to the unique solution of \eqref{eq:probuniform-problem}  with $p=1$. If the objective happens to be even strongly convex (as in the numerical examples), then any bounded sequence $(u_{p_k})_{k\in \mathbb{N}}$ of minimizers of problem \eqref{eq:probuniform-problem}  with $p_k \in (0,1)$ strongly converges to the unique solution of problem \eqref{eq:probuniform-problem}  with $p=1$.
\end{proposition}
\proof
As proven in Proposition \ref{convexity}, problem \eqref{eq:probuniform-problem}  has a convex feasible set for all $p \in (0,1]$. Therefore, if $F$ is strictly convex, the problem has a unique minimizer. Hence, if $u_{p_k}$ is bounded, then as per Lemma \ref{prep}, the only possible accumulation point is the unique minimizer of \eqref{eq:probuniform-problem} for $p=1$. 

Strong convexity of $F$ (with modulus $\mu$) gives 
\begin{equation}
\label{eq:strong-convexity}
F(u_{p_k}) - F(\bar{u}) \geq\langle\nabla F(\bar{u}), u_{p_k}-\bar{u}\rangle + \frac{\mu}{2}\lVert u_{p_k}-\bar{u}\rVert^2.
\end{equation}
From Lemma \ref{prep}, we have 
$\lim_{k\rightarrow \infty} F(u_{p_k})= F(\bar{u})$. Also $\langle\nabla F(\bar{u}), u_{p_k}-\bar{u}\rangle\to 0$ as $k\rightarrow \infty$ by weak convergence of $u_{p_k}$. Strong convergence of the same sequence follows by \eqref{eq:strong-convexity}.
\qed

\bigskip\noindent
Note that the last Proposition explains why in Figure \ref{fig1a} we observe strong and not just weak convergence. 
\subsection{Moreau--Yosida approximation for dealing with almost sure constraints}
A difficulty in the numerical solution to problem \eqref{eq:a.s.-problem} is in the enforcement of the state constraints. One possibility to handle this computationally is to penalize this constraint in the objective. Note that $(H_0^1(D), L^2(D), H^{-1}(D))$ is a Gelfand triple. Moreover, the cone $\mathcal{K}= \{ y \in L_{\pP}^\infty(\Omega,H_0^1(D)) \mid y(x,\omega) \leq 0 \text{ a.e.-$\pP$-a.s.}\}$ is compatible with the cone $\mathcal{K}_H = \{ y \in L_{\pP}^2(\Omega,L^2(D)) \mid y(x,\omega) \leq 0 \text{ a.e.-$\pP$-a.s.}\}$ in the sense that $\mathcal{K}_H \cap L^\infty_{\pP}(\Omega,H_0^1(D)) = \mathcal{K}$. This justifies the penalization on the weaker (Hilbert) space $H=L_{\pP}^2(\Omega,L^2(D))$ where the projection $\pi_{\mathcal{K}_H}$ onto the cone $\mathcal{K}_H$ is possibly computationally cheaper. The Moreau--Yosida regularization (or envelope) for the indicator function on the cone $\mathcal{K}_H$ has the formula
\[
\hat{\beta}^\gamma(k) = \gamma\lVert k - \pi_{\mathcal{K}_H}(k)\rVert_{H}^2;
\]
see, e.g., \cite[Section 4.1]{Geiersbach2022a}. Let $\boldsymbol{\alpha} \in L^\infty_{\pP}(\Omega, H^1(D))$ be the function that is equal to $\alpha$ a.e.~in $D$ and $\pP$-a.s. Penalizing the constraint \eqref{eq:a.s.-problem-constraint3} amounts to adding $\hat{\beta}^\gamma(y-\boldsymbol{\alpha})$ to the objective, leading to the modified problem 
\begin{subequations}
    \label{eq:a.s.-problem-MY}
    \begin{alignat}{3}
    \min_{u \in L^2(D)}\,  \Big\lbrace f^\gamma(u) :=F(u)+\gamma \E[\lVert &\max(0,y-\boldsymbol{\alpha})\rVert^2_{L^2(D)}] \Big\rbrace & &&  \\
   \text{s.t.}  \quad -\Delta  y(x,\omega)  &=  u(x) + f(x,\xi (\omega)), & &&\quad x \in D \quad \text{$\mathbb{P}$-a.s.}, \label{eq:a.s.-problem-constraint1-MY} \\
 y(x,\omega) &=0, & &&\quad x \in \partial D  \quad \text{$\mathbb{P}$-a.s.} \label{eq:a.s.-problem-constraint-bc}
    \end{alignat}
\end{subequations}
The results from \cite{Geiersbach2022a} focus on the consistency of the optimality conditions for \eqref{eq:a.s.-problem-MY} to the optimality conditions for \eqref{eq:a.s.-problem}, which require an interior point condition (and therefore higher regularity of the solution to the PDE). Here, we will focus on the consistency of the primal problem and will work with the weakest regularity available. 
Let $A \colon H_0^1(D) \rightarrow H^{-1}(D)$ represent the Laplacian,  which is a linear isomorphism thanks to Assumption~\ref{ass:PDE-standing} and the Lax--Milgram lemma. We denote with $\mathcal{A}^{-1}\colon L_{\pP}^\infty(\Omega,L^2(D)) \rightarrow L_{\pP}^\infty(\Omega, H_0^1(D))$ the superposition defined by $[\mathcal{A}^{-1}y](\omega)=A^{-1}y(\cdot,\omega).$ Let $\mathcal{B}\colon L^2(D) \rightarrow L_{\pP}^\infty(\Omega,L^2(D))$ be the canonical embedding. Let $\tilde{f}$ be defined by $\tilde{f}(x,\omega):=f(x,\xi(\omega)).$
It will be convenient to define the (affine linear) control-to-state operator $\tilde{\mathcal{S}}\colon L^2(D) \rightarrow L_{\pP}^\infty(\Omega, H_0^1(D))$ by $\tilde{\mathcal{S}}(u) = \mathcal{A}^{-1}(\mathcal{B}u + \tilde{f})$, which is bounded by \cite[Lemma 10]{geiersbach2023optimality}. Due to the continuous canonical embedding $\iota \colon L^\infty_{\pP}(\Omega,H_0^1(D)) \rightarrow L^2_{\pP}(\Omega,L^2(D))$, the operator $\mathcal{S}:=\iota \circ \tilde{\mathcal{S}}$ is continuous from $L^2(D)$ to $L^2_{\pP}(\Omega,L^2(D)).$

We demonstrate how solutions to \eqref{eq:a.s.-problem-MY} converge to \eqref{eq:a.s.-problem} in the limit as $\gamma \rightarrow \infty$. 
First we need the following.
\begin{lemma}
\label{lemma:continuity properties}
For any $\gamma\geq 0$, the reduced functional $\beta^\gamma \colon L^2(D) \rightarrow \R$ defined by
\[
\beta^\gamma(u) =  \gamma\E[\lVert \max(0,\mathcal{S}(u)-\boldsymbol{\alpha})\rVert_{L^2(D)}^2]
\]
is convex and (weakly) lower semicontinuous.
\end{lemma}
\proof
Since $\mathcal{S}$ is continuous and the map $y \mapsto \E[\lVert \max(0,y-\boldsymbol{\alpha})\rVert_{L^2(D)}^2]$ is evidently continuous on $L^2_{\pP}(\Omega, L^2(D))$, $\beta^\gamma$ is continuous on $L^2(D)$. Moreover, since $\mathcal{S}$ is (affine) linear, $\beta^\gamma$ is convex. Weak lower semincontinuity follows from the continuity and convexity of $\beta^\gamma$.
\qed

\bigskip\noindent
We now have the following result: 
\begin{lemma}
\label{lem:MY-consistency}
If $F$ is (in addition to the usual assumptions) coercive, then there exists a solution $u^\gamma$ to \eqref{eq:a.s.-problem-MY} for all $\gamma >0.$ Additionally, if there exists a feasible point for problem \eqref{eq:a.s.-problem}, then given a sequence $( \gamma_n)$ with $\gamma_n \rightarrow \infty$, weak limit points of $u^{\gamma_n}$ solve \eqref{eq:a.s.-problem}. 
\end{lemma}
\proof
Let $F_{\textup{ad}}= \{ u \in L^2(D) \mid [\mathcal{S}(u)](x,\omega) \leq \alpha \text{ a.e.-$\pP$-a.s.}\}$ denote the feasible set for problem \eqref{eq:a.s.-problem}. It is straightforward to show that this set is convex; it is closed thanks to the continuity of $\mathcal{S}.$ Thus, $F_{\textup{ad}}$ is weakly closed (see \cite[Theorem 2.23]{Bonnans2013}). Coercivity of $F$ (in combination with the other assumptions) implies that $f^\gamma$ attains its minimum over $F_{\textup{ad}}$. 
In particular, there exists a minimizing sequence $(u_{n})_{n \in \mathbb{N}} \subset F_{\textup{ad}}$ such that $\lim_{n \rightarrow \infty} f^\gamma(u_n) = \inf_{u \in F_{\textup{ad}}} f^\gamma(u).$ This sequence is also bounded due to the coercivity of $F$, so we can extract a subsequence $(u_{n_k})_{k \in \mathbb{N}}$ such that $u_{n_k} \rightharpoonup u^\gamma$. Now, we have (since $F$ is weakly lower semicontinuous, as convex and continuous function)
\begin{align*}
    \liminf_{k \rightarrow \infty} F(u_{n_k}) + \beta^\gamma(u_{n_k}) \geq     \liminf_{k \rightarrow \infty} F(u_{n_k}) + \liminf_{k \rightarrow \infty} \beta^\gamma(u_{n_k}) \geq F(u^\gamma)+\beta^\gamma(u^\gamma),
\end{align*}
meaning $u^\gamma$ solves problem \eqref{eq:a.s.-problem-MY}.
In other words, a solution exists for any $\gamma>0$.

Now, let $(u^{\gamma_n})_{n \in \mathbb{N}}$ be a sequence such that $\gamma_n \rightarrow \infty$ as $n \rightarrow \infty$. Let $\bar{u} \in F_{\textup{ad}}\neq \emptyset$ be a minimizer to \eqref{eq:a.s.-problem}. Then, by optimality of $u^{\gamma_n}$ and feasibility of $\bar{u}$, we have
\begin{equation}
    \label{eq:MY-inequalities-limit-argument}
    F(u^{\gamma_n}) \leq F(u^{\gamma_n}) +\beta^{\gamma_n}(u^{\gamma_n}) \leq F(\bar{u}) + \beta^{\gamma_n}(\bar{u}) = F(\bar{u}).
\end{equation}
This means $u^{\gamma_n}$ is a minimizing sequence for $F$, which is bounded thanks to the coercivity of $F$. Therefore, there exists a subsequence $(u^{\gamma_{n_k}})_{k \in \mathbb{N}}$ and point $\hat{u}$ such that $u^{\gamma_{n_k}} \rightharpoonup \hat{u}.$ Weak lower semicontinuity of $F$ implies that
\begin{equation}
\label{eq:weak-lsc-inequality-MY}
F(\hat{u}) \leq \liminf_{k \rightarrow \infty} F(u^{\gamma_{n_k}}) \leq F(\bar{u}).
\end{equation}
It only remains to show that $\hat{u}$ is feasible for \eqref{eq:a.s.-problem}. Since $F$ is bounded over $F_{\textup{ad}}$, \eqref{eq:MY-inequalities-limit-argument} implies that there exists a constant $c>0$ such that
\[
\beta^{\gamma_{n_k}}(u^{\gamma_{n_k}}) \leq c \quad \Leftrightarrow \quad \E[\lVert \max(0,\mathcal{S}(u^{\gamma_{n_k}})-\boldsymbol{\alpha})\rVert_{L^2(D)}^2] \leq \frac{c}{\gamma_{n_k}}. 
\]
Since $\gamma_{n_k} \rightarrow \infty$, we obtain $\E[\lVert \max(0,\mathcal{S}(u^{\gamma_{n_k}})-\boldsymbol{\alpha})\rVert_{L^2(D)}^2] \rightarrow 0$, implying that $[\mathcal{S}(\hat{u})](x,\omega)\leq \alpha$ a.e.~$\pP$-a.s. By definition of $\mathcal{S}$, the function $\mathcal{S}(u)$ satisfies \eqref{eq:a.s.-problem-constraint1}--\eqref{eq:a.s.-problem-constraint2}.

\qed

\begin{corollary}
Under the same conditions as Lemma~\ref{lem:MY-consistency}, suppose that $F$ is strongly convex. Then any bounded sequence $(u^{\gamma_n})_{n\in \mathbb{N}}$ of minimizers of problem \eqref{eq:a.s.-problem-MY} with $\gamma_n \rightarrow \infty$ converges to the unique solution of problem \eqref{eq:a.s.-problem}.
\end{corollary}

\proof
Since there can only be one solution to \eqref{eq:a.s.-problem}, \eqref{eq:weak-lsc-inequality-MY} together with the feasibility of $\hat{u}$ and $\bar{u}$ imply that $\lim_{n \rightarrow \infty} F(u^{\gamma_n}) = F(\bar{u})$.
Strong convergence follows with the same reasoning used in Proposition~\ref{propoPartial}.
\qed

\begin{remark}
While the above proof is based on \cite[Proposition 3.8]{Geiersbach2022a}, we used the reduced formulation here, which greatly simplified certain arguments; we did not need to rely on arguments using the weak* topology, for one. Moreover, we did not require an interior point condition to establish consistency (which was assumed in \cite{Geiersbach2022a} but not used in this proof). 
\end{remark}
\subsection{Numerical results by ``sampling the distribution,'' ``sampling the support,'' or ``sampling the boundary'' of the support}\label{saamy}
\begin{figure}[htb]
    \centering
    \includegraphics[width=7cm]{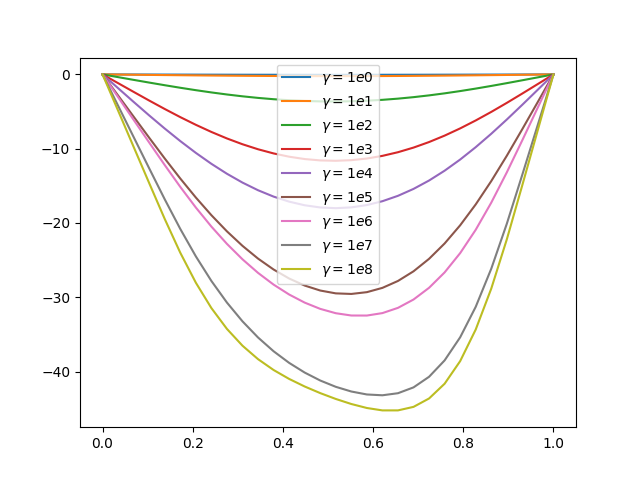}
    \hspace{-0.8cm}
    \includegraphics[width=7cm]{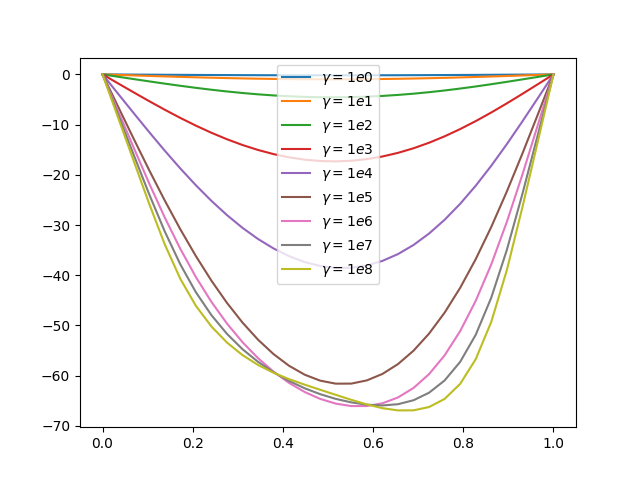}
    
    \includegraphics[width=7cm]{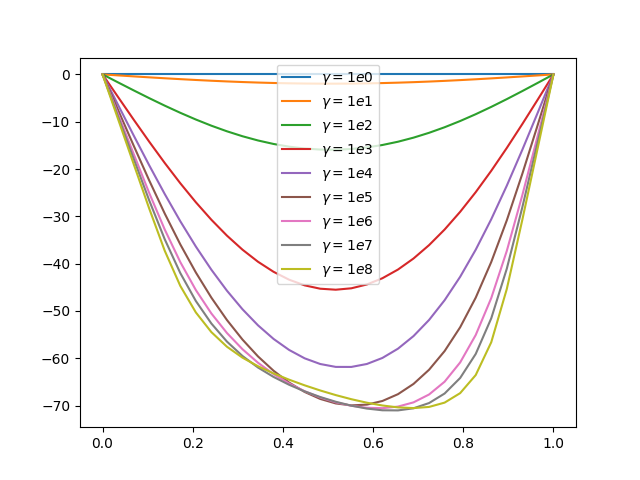}
    \hspace{-0.4cm}
    \includegraphics[width=6.4cm]{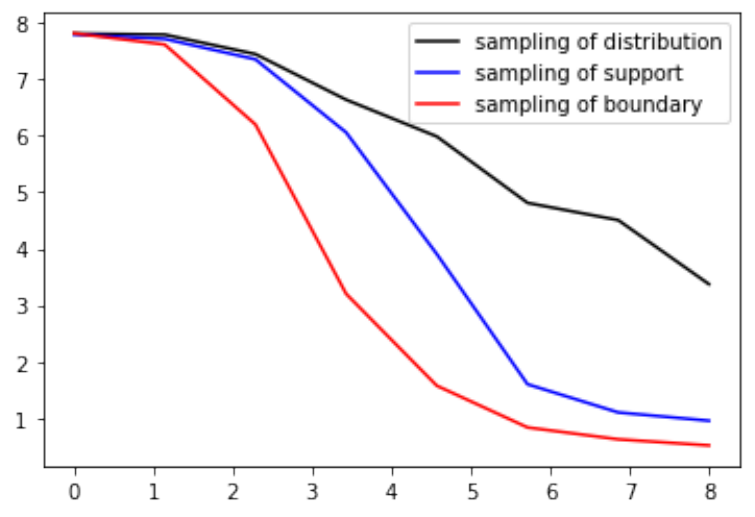}
    \caption{Solutions $u^\gamma$ obtained for increasing values $\gamma =10^k$ and increasing sample sizes $3^k$ using sampling of the distribution (top left),  sampling of the support (top right), and sampling on the boundary of the support (bottom left). The decrease of constraint violation as a function of $k$ for the three sampling approaches is illustrated in the diagram at the bottom right.}
    \label{figMY1D}
\end{figure}

In the following numerical tests, we use gradient descent to solve a sample average approximation of the subproblem  \eqref{eq:a.s.-problem-MY}. Given a fixed number $N$ of samples $z_1, \dots, z_{N}$ of the random vector $\xi$, the corresponding SAA problem is
\begin{subequations}
    \label{eq:a.s.-problem-MY-SAA}
    \begin{align}
    \min_{u \in L^2(D)}\,  \Big\lbrace \hat{f}^{\gamma}(u) :=F(u)+\frac{\gamma}{N} \sum_{i=1}^{N}[\lVert \max(0,y_i-\boldsymbol{\alpha})\rVert^2_{L^2(D)}]  \Big\rbrace&    \\
   \text{s.t.}  \quad -\Delta  y_i(x)  =  u(x) + f(x,z_i),\quad x \in D \quad \forall i=1, \dots, N,& \label{eq:a.s.-problem-constraint1-MY-SAA} \\
 y_i(x) =0, \quad x \in \partial D  \quad \forall i=1, \dots, N.& \label{eq:a.s.-problem-constraint-bc-SAA}
    \end{align}
\end{subequations}
The gradient of $\hat{f}^{\gamma}$ can be computed using standard techniques, resulting in 
\begin{equation}
    \nabla \hat{f}^{\gamma}(u) = \nabla F(u) - \frac{1}{N} \sum_{i=1}^{N} p_i,
\end{equation}
where $p_i$ is the solution to the adjoint equation
\begin{equation}
\label{eq:adjoint}
    -\Delta p_i(x) = -2 \gamma
\max (0, y_i(x)-\alpha), \quad x \in D \quad \text{and} \quad p_i(x) = 0 \quad x \in\partial D,
\end{equation}
and $y_i$ is the solution to \eqref{eq:a.s.-problem-constraint1-MY-SAA}--\eqref{eq:a.s.-problem-constraint-bc-SAA} for a fixed sample $z_i$. 
\begin{algorithm}
\begin{algorithmic}[0]
\caption{SAA Gradient descent with Moreau--Yosida penalization}\label{alg} 
\STATE {\bf Input}: $D\subset\mathbb{R}^d$ (domain), $K$ (number of penalty updates), $\{z_{1}, \dots, z_{N_{K}}\}$ (set of samples), $u_1$ (initial control), $(\gamma_k)_{k\in \{1, \dots, K \}}$ (sequence of penalties), $(t_j)_j$ (step-size rule)
\FOR{$k=1, \dots, K$}
\FOR{$j=1, 2, \dots $}
        \IF{$\lVert \nabla \hat{f}^{\gamma_k}(u_j)\rVert > \textup{tol}$}
        \STATE break
        \ENDIF
        \FOR{$\ell =1, \dots, N_k$}
        \STATE $y_{\ell} \gets$ solution to \eqref{eq:a.s.-problem-constraint1-MY-SAA}--\eqref{eq:a.s.-problem-constraint-bc-SAA} with $z_{\ell}$ and $u_j$ 
        \STATE $p_\ell \gets$ solution to \eqref{eq:adjoint} with $y_{\ell}$ and $\gamma_k$
        \ENDFOR
        \STATE $u_{j+1}:=u_j - t_j (\nabla F(u_j) - \frac{1}{N_k} \sum_{i=1}^{N_k} p_i)$
\ENDFOR
\STATE $u_1:=u_{j}$
\ENDFOR
\STATE {\bf Output:} Solution $u$.
\end{algorithmic}
\end{algorithm}
For numerical tests, we use the same setup as described at the beginning of Section~\ref{Sect_resulstindimone}. However, rather than considering the purely Gaussian distribution $\mathcal{N}(0,\Sigma)$ defined there, we truncate it to the ellipsoid $\mathcal{E}$ defined below the original data in connection with Fig. \ref{fig1} (right). The sampling of this truncated Gaussian distribution is carried out by accepting/rejecting samples of the underlying original Gaussian distribution according to whether or not these samples belong to the ellipsoid $\mathcal{E}$.
The PDEs \eqref{eq:a.s.-problem-constraint1-MY-SAA}-- \eqref{eq:a.s.-problem-constraint-bc-SAA} and \eqref{eq:adjoint} were solved using FEniCS \cite{Alnes2015} with a finite element discretization over 29 intervals. 
We use a path-following approach, where in an outer iteration, the penalty was increased and in inner iterations, gradient descent on a subproblem was performed until a given tolerance. A starting value for the control was set to $u\equiv -1$. In each outer iteration $k \in \{ 0, \dots, 8\}$, $\gamma = 10^k$, an increasing batch of $3^k$ $m$-dimensional vectors was used. The step-size for each outer iteration was chosen to be $t_\ell = 4/\ell$, which was informed by the strong convexity of the problem; see \cite{Geiersbach2019a} for a discussion of this choice in the context of PDE-constrained optimization. Each inner iteration was terminated when $\lVert \nabla \hat{f}^\gamma(u_j)\rVert_{L^2(D)} < 10^{-4}$. 
In Fig. \ref{figMY1D} (top left), different optimal controls are displayed for increasing values $\gamma=10^k$ and increasing sample sizes $3^k$. 

So far, we have considered samples of the distribution
in order to characterize solutions of the almost sure model. Sometimes, however, one can do better than sampling the distribution. As observed in \cite[Example 3.3]{geiersbach2023optimality}, problems with almost sure constraints are not equivalent, in general, to problems with worst-case constraints on the support of the random parameter \cite[Lemma 3.4]{geiersbach2023optimality}. They are, however, if the underlying random inequality is lower semicontinuous with respect to the random parameter. This is the case, for instance, in our problem. Therefore, the problem \eqref{eq:a.s.-problem} with almost sure constraints is equivalent with the \textbf{robust} optimization problem
\begin{subequations}
    \label{eq:robust-problem}
    \begin{alignat}{3}
    \min_{u \in L^2(D)}\, F(u) & & &&  \\
   \text{s.t.}  \quad -\Delta  \hat{y}(x,z)  &=  u(x) + f(x,z), & &&\quad (x,z) \in D \times \Xi, \label{eq:uniform-problem-constraint1} \\
 \hat{y}(x,z) &=0, & &&\quad (x,z) \in \partial D \times \Xi, \label{eq:uniform-problem-constraint2}\\
 \hat{y}(x,z), &\leq \alpha & && (x,z)\in D \times \Xi, \label{eq:uniform-problem-constraint3}
    \end{alignat}
\end{subequations}
where the support $\Xi$ of the random vector $\xi$ plays the role of the uncertainty set with respect to which the worst case has to be taken into account. Since the support $\Xi$---contrary to our concrete simple ellipsoid $\mathcal{E}$---can be potentially complicated and hard to deal with, one approach might consist in {\it uniformly} sampling this support. Of course, our previous sampling of the distribution trivially provides samples on the support. These are, however, strongly concentrated around the mean $0$ of the ellipsoid due to the underlying Gaussian distribution, while only few of these samples will be close to the boundary of the ellipsoid. However, points at the boundary are much more informative in defining the robust constraint (see discussion in Section \ref{samplefree} below).
Therefore, it seems reasonable, if possible, to work with a uniform sample of the support or even just its boundary. In the case of our ellipsoid $\mathcal{E}:=\{z\in\mathbb{R}^m\mid z^T\Sigma^{-1}z\leq 36\}$, a uniformly distributed sample can be created as follows: Let $\{v^{(k)}\}_{k=1}^K$ be a sample of the uniform distribution on the sphere $\mathbb{S}^{m-1}$ as in Section \ref{spherapprox} and let $\{\tau^{(k)}\}_{k=1}^K$ be a sample of the uniform distribution on $[0,1]$. Then, $\{6\tau^{(k)}\Sigma^{1/2}v^{(k)}\}_{k=1}^K$ is a sample of the uniform distribution on $\mathcal{E}$. When fixing $\tau^{(k)}\equiv 1$ instead, one rather obtains a uniform distribution on the boundary of $\mathcal{E}$.
With this potentially more efficient uniform sampling of the support $\Xi$ or its boundary, respectively, one may numerically proceed in the same way as with the previous sampling of the distribution. The corresponding results are illustrated in Fig. \ref{figMY1D} (top right and bottom left). A plot showing the decrease in constraint violation for the three sampling strategies is illustrated in the diagram at the bottom right. More precisely, constraint violation refers to the maximum excess over the threshold $\alpha$ (with respect to the domain of the space variable $x$ and to the support of the random variable $\xi$) of the random state under optimal control.
All sampling schemes show convergence according to the theoretical results, but uniform sampling of the support and even more of just its boundary seem to exhibit much faster convergence than sampling of the distribution (see amplitudes of solutions in Fig. \ref{figMY1D}). This will be confirmed in the following section upon comparing results with the sample-free solution of the problem.
\subsection{Sample-free solution of the optimization problem with almost-sure constraints}\label{samplefree}
So far, we have seen four alternative approaches for dealing with almost sure constraints: one by formulating a probabilistic constraint with probability level $p=1$ and three sample average approaches based on Moreau--Yosida approximation with increasing penalty parameter. The latter three methods differed according to whether the distribution was sampled itself or rather its support or even just the boundary of the support. In this section, we present yet another alternative, namely the direct solution of the robust optimization problem \eqref{eq:robust-problem} without sampling. This solution can be understood as the ``true'' solution of the almost sure problem. We emphasize that such a sample-free approach may not be possible in general, when the objective and the support are not simple enough to provide an analytical solution of the inner problem below.

We recall our assumption that the support $\Xi\subseteq\mathbb{R}^m$ of the random vector $\xi$ is compact.
Similar to the probabilistic setting, we may exploit the parametric control-to-state operator $S$ there, in order to recast \eqref{eq:robust-problem} as an infinite-dimensional optimization problem
\begin{equation}\label{optprob2}
\min\limits_{u\in L^2(D)}\,\,F(u)\mbox{ s.t. }h(u)\leq 0,\quad h(u):=\max\limits_{(x,z)\in\bar{D}\times\Xi }[S(u,z)](x)-\alpha.
\end{equation}
 In order to solve \eqref{optprob2} numerically, we need to compute the function $h$ and its (sub-) gradients. It is advantageous to rewrite $h$ as an iterated maximum
\begin{equation}\label{doublemax}
h(u)=\max\limits_{x\in\bar{D}}\max\limits_{z\in\Xi }[S(u,z)](x)-\alpha.
\end{equation}
Since $S$ is affine linear in $z$, the inner maximization problem (for some fixed $u\in L^2(D)$ and $x\in\bar{D}$) can be solved analytically if the support $\Xi$ is a simple set (e.g., rectangle, ellipsoid). Then, since the domain $D$ is of low dimension, the outer maximization is easily carried out on a grid corresponding to the discretization of the state. For instance, if the support is given by an ellipsoid 
\begin{equation}\label{ellipsoid}
\Xi =\{z\in\mathbb{R}^m\mid z^TBz\leq\gamma\}
\end{equation}
for some symmetric and positive definite matrix $B$ and some $\gamma >0$, then a linear function $c^Tz$ realizes its maximum over $\Xi$ at the point 
\[
z^*=\sqrt{\frac{\gamma}{c^TB^{-1}c}}B^{-1}c.
\]
Now, in order to compare the sample-free robust solution with the previous three methodologies, assume that the support of $\Xi$ of the truncated Gaussian random vector $\xi$ is the ellipsoid $\mathcal{E}:=\{z\in\mathbb{R}^m\mid z^T\Sigma^{-1}z\leq 36\}$ as in Section \ref{Sect_resulstindimone}. 
For a fixed $u\in L^2(D)$, denote as in Theorem \ref{exactclarke} by $\bar{y}_u:=S(u,0)$ the mean state associated with $u$, i.e., the solution of the PDE \eqref{eq:probuniform-problem-c}, \eqref{eq:probuniform-problem-d} with right-hand side $u(x)+f_0(x)$. Then, with \eqref{sstructure}, we infer that
\[
[S(u,z)](x)=\bar{y}_{(u)}(x)+\sum\limits_{i=1}^{m}z_iy^{(i)}(x)\quad\forall(x,z)\in\bar{D}\times\Xi .
\]
Hence, for some fixed $u\in L^2(D)$ and $x\in\bar{D}$, the inner maximization problem mentioned above amounts to maximizing the linear function $c^Tz$ with $c:=(y^{(i)}(x))_{i=1}^m$ over the ellipsoid \eqref{ellipsoid} with $B:=\Sigma^{-1}$ and $\gamma:=36$, where the data for the support have been chosen as in the previous computations. Therefore, the maximum is realized for
\[
z^*=\frac{6}{\sqrt{c^T\Sigma c}}\Sigma c.
\]
Summarizing, the function value of $h$ can be determined as 
\[
h(u)=\,\max\limits_{x\in\bar{D}}\,\bar{y}_{(u)}(x)+\sum\limits_{i=1}^{m}z^*_iy^{(i)}(x)-\alpha =:\max\limits_{x\in\bar{D}}H(u,x)-\alpha,
\]
which is easily approximated over the discretized state in low dimensions. The robust constraint function $h$ will be typically nonsmooth (see discussion below). Its (convex) subdifferential is easily characterized as
\[
\partial h(u)={\rm clco}\,\{{\bf u}_x\mid x\in M(u)\};\quad M(u):=\{x\in\bar{D}\mid h(u)=H(u,x)-\alpha\} ,
\]
where ${\bf u}_x$ is introduced below \eqref{basicstates}. 
Equipped with these tools, an appropriate algorithm for nonsmooth convex optimization can be employed in order to numerically solve problem \eqref{optprob2}.
\begin{figure}[htb]
    \centering
    \includegraphics[width=5.8cm]{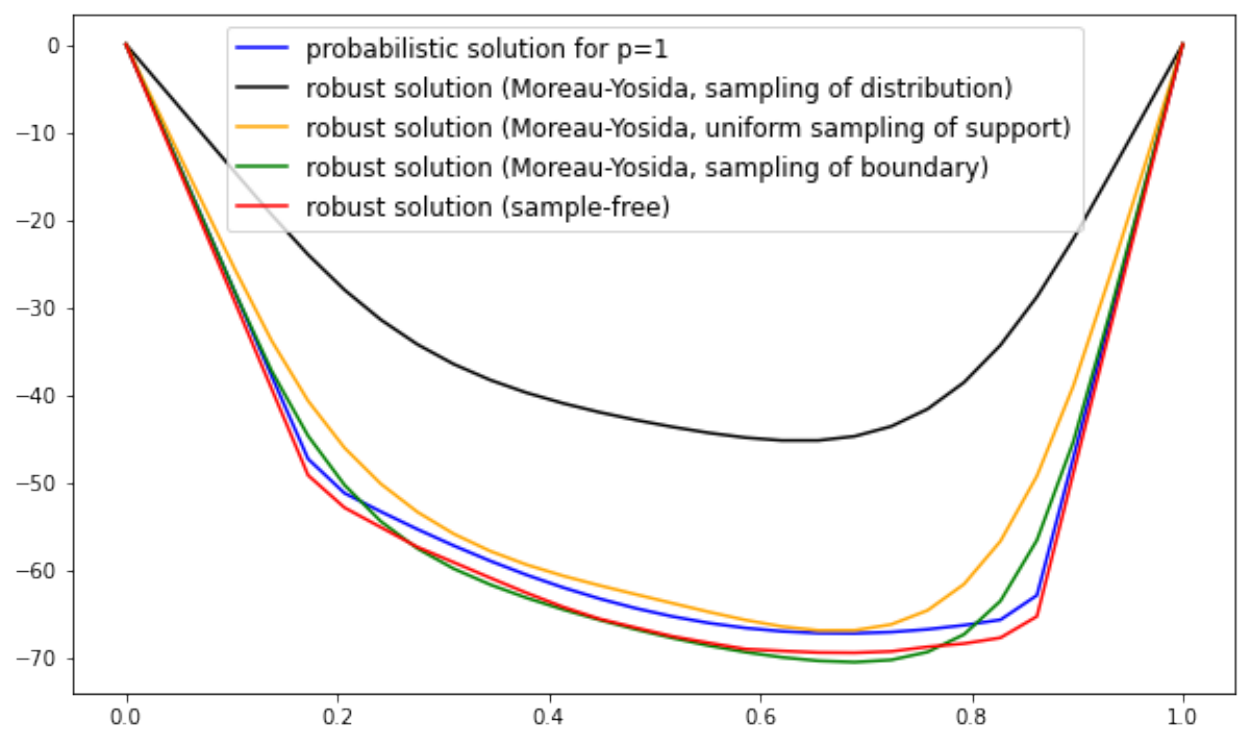}
    \includegraphics[width=5.8cm]{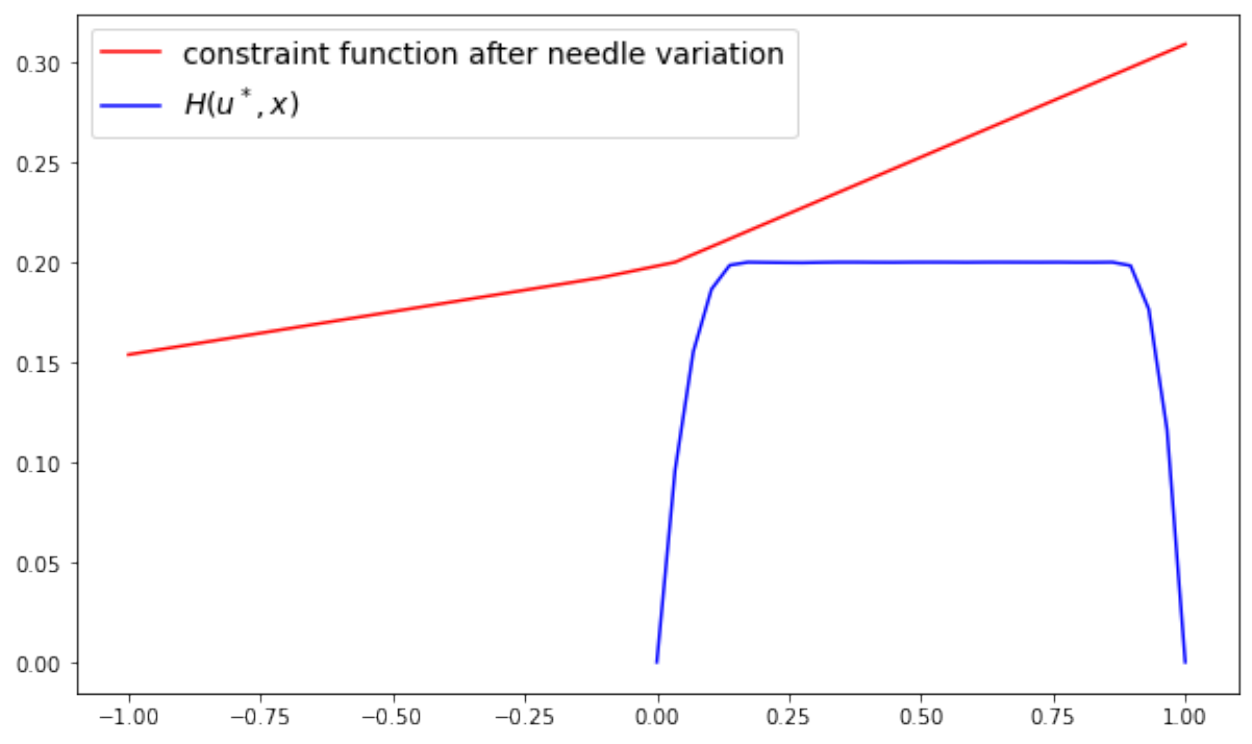}
    \caption{Comparison of the probabilistic, three sampling-based the sample-free  solutions to the worst-case problem (left). Illustration of the nonsmoothness of the robust constraint function $h$ (right), details see text.}
    \label{fig2}
\end{figure}
Figure \ref{fig2} (left) shows the comparison of the probabilistic solution in Figure \ref{fig1} (right) for $p=1$ with the three sampled and the sample-free robust solutions on the same support. The sampled robust solutions are those illustrated in Fig. \ref{figMY1D} with penalty parameter $\gamma=10^8$ and sample size $3^8$. Several conclusions can be drawn when understanding, as suggested above, the sample-free solution as the ``true'' reference: First, sampling the distribution is the least efficient method. This is not surprising, because in this sampling procedure a minimum amount of information is exploited. In particular, no knowledge about the support of the distribution is used. Sampling the support uniformly (and not according to the original random distribution, which might be strongly concentrated in the interior of the support) yields a much better approximation. This can even be significantly improved when just sampling uniformly the boundary of the support. 
The explanation of this further improvement is that, thanks to the linearity of the objective in the inner maximization problem \eqref{doublemax}, all constraints are already induced by boundary points. This observation from Fig. \ref{fig2} complements those from Fig. \ref{figMY1D}.
Best performance in approximating the ``true'' solution is reached by the just-mentioned sampling of the boundary of the support combined with a Moreau--Yosida approximation on the one side and probabilistic constraint with level $p=1$ on the other side. The difference between both results is that the former method better approximates the smooth pieces of the solution, whereas the latter one is more precise close to the nonsmooth kinks.
\begin{remark}
We note that the benefit of sampling the boundary of the support pertains to constraint functions $g(u,z)$ that are not necessarily linear but just convex in the random variable $z$. Then, the objective in the inner maximization of \eqref{doublemax} is convex and, hence, assuming that the support $\Xi$ is not just compact but also convex, the maximizers are still located on the boundary of $\Xi$.
\end{remark}
\noindent
It is interesting to observe that the robust constraint function $h$ fails to be differentiable at the optimal solution $u^*$ (here, we do not refer to the fact, evident from Fig. \ref{fig2} (left), that this solution itself is a nonsmooth function). A first indication of this is that the function $H(u^*,\cdot)$ plotted in Fig. \ref{fig2} (right) achieves its maximum over $\bar{D}=[0,1]$ on a whole interval. Hence, the set $M(u^*)$ and the subdifferential $\partial h(u^*)$ do not shrink to a 
singleton and, consequently, $h$ is nonsmooth. This fact can be alternatively illustrated by plotting $h$ as a univariate function of a needle variation of $u^*$ on a small interval $\Delta\subseteq D$ in the middle of the domain. More precisely, we define the function
$\tilde{h}(t):=h(u_t)$, where 
\[
u_t(x):=\left\{\begin{array}{ll}
     u^*(x)+t  &\mbox{ if }x\in\Delta\\
     u^*(x)&\mbox{ if } x\in D\setminus\Delta
\end{array} \right.
\]
Evidently, if $h$ was differentiable, then so would the univariate function $\tilde{h}$ be, too. This, however, is not the case as can be seen from
Fig. \ref{fig2} (right). Note that in the figure, we are simultaneously using the horizontal axis for the $x$-variable of the function $H(u^*,x)$ and for the $t$-variable of the function $\tilde{h}$.
\section*{Acknowledgments}

The first two authors thank the Deutsche Forschungsgemeinschaft for their support within projects B02 and B04 in the “Sonderforschungsbe\-reich/Transregio 154 Mathematical Modelling, Simulation and Optimization using the Example of Gas Networks.” The first author acknowledges support from the Deutsche Forschungsgemeinschaft (DFG, German Research Foundation) under Germany's Excellence Strategy – The Berlin Mathematics Research Center MATH+ (EXC-2046/1, project ID: 390685689). The second
author acknowledges support by the FMJH Program Gaspard Monge in optimization and operations research including support to this program by EDF. The third author was supported by  Centro de Modelamiento Matem\'atico (CMM), ACE210010 and FB210005, BASAL funds for center of excellence and ANID-Chile grant: Fondecyt Regular 1200283 and Fondecyt Regular 1220886 and Proyecto Exploración 13220097.
\appendix
\section{Appendix}
In the following, we provide a proof of Theorem \ref{clarkeequality}. It will follow from two lemmas below.
\begin{lemma}
Under the assumptions of Theorem \ref{clarkeequality}, the mapping $g$ defined there, satisfies the assumptions of Theorem \ref{probderiv}.
\end{lemma}
\proof 
We refer to the assumptions 1.-4. of Theorem \ref{clarkeequality}. The continuity of $h$ (see 1.) ensures that $g$ can be written as a max rather than just sup and that $g(\bar{u},0)<0$ by 3. Condition 2. implies that $g$ is locally Lipschitzian. Indeed, the function $\gamma (u,z):=\max_{w\in K}\|D_{(u,z)}h(u,z,w)\|$ is continuous by 2., hence it is locally bounded around an arbitrary $(u,z)$ by some constant $C$. This implies that all functions $h(\cdot ,\cdot ,w)$ for $w\in K$ are locally  Lipschitz continuous around $(u,z)$ with a common Lipschitz constant $C$. Hence, $g$ is locally Lipschitz continuous around $(u,z)$. The convexity of $h$ in its second argument (see 1.) yields that $g$ is convex in its second argument too. Therefore, our general assumption (GA) on $g$ holds true. We show that 4. implies \eqref{growth}. To this aim, for $(u,d,z,t)\in U\times U\times\mathbb{R}^m\times\mathbb{\R}$, denote by 
$w_{u,d,z,t}\in K$ some element with $g(u+td,z)=h(u+td,z,w_{u,d,z,t})$, which exists thanks to the compactness of $K$ and the continuity of $h$. It follows that $g(u'+td,z)=h(u'+td,z,w_{u',d,z,t})$ and $-g(u',z) \leq - h(u',z,w_{u',d,z,t})$, which yields 
\begin{eqnarray*}
g^{\circ}(\cdot ,z)(u;d)&=&\limsup\limits_{u'\to u,t\downarrow 0}\frac{g(u'+td,z)-g(u',z)}{t}\\&\leq&\limsup\limits_{u'\to u,t\downarrow 0}\frac{h(u'+td,z,w_{u',d,z,t})-h(u',z,w_{u',d,z,t})}{t}\\
&\leq&\limsup\limits_{u'' \to u}\sup\limits_{w\in K}\|D_uh(u'',z,w)\| \| d\|\\&\leq&c\exp({\|z\|})\|d\|\quad\forall u:\|u-\bar{u}\|\leq (2c)^{-1}\quad\forall z:\|z\|\geq c\quad\forall d\in U. 
\end{eqnarray*}
Since, for some $c_0$,
\[
\exp\left({\|z\|}\right)\leq\left\Vert z\right\Vert^{-m}\exp\left(\frac{\left\Vert z\right\Vert^2}{2\left\Vert \Sigma^{1/2}\right\Vert^2} \right)\quad\forall z:\|z\|\geq c_0,
\]
\eqref{growth} holds true with $l:=\max\{1, c_0,2c\}$.

\qed

\bigskip\noindent
As a consequence of the previous lemma, Theorem \ref{clarkeequality} yields \eqref{clarkeinclu} via Theorem \ref{probderiv}. The argument provided in \cite{Hantoute2019} to infer \eqref{clarkeinclu} from the assumptions of Theorem \ref{probderiv}, relies on the functions $e(\cdot ,v)$ being uniformly (with respect to all $v\in\mathbb{S}^{m-1}$) Lipschitzian with a common modulus on a neighborhood of $\bar{u}$, see \cite[Theorem 5, Corollary 2]{Hantoute2019}. Since constants are integrable with respect to the uniform distribution on the sphere, one may invoke Clarke's theorem on subdifferentiation of integral functionals \cite[Theorem 2.7.2]{clarke1983} in order to derive 
\eqref{clarkeinclu}. According to an addendum in Clarke's theorem, equality in \eqref{clarkeinclu} along with Clarke regularity of $\varphi$ at $\bar{u}$ could be inferred if the partial functions $e(\cdot ,v)$ were Clarke regular at $\bar{u}$. Unfortunately, they are not. On the other hand, with $e(\cdot ,v)$ also the negative functions  $-e(\cdot ,v)$ are uniformly Lipschitzian with a common modulus on a neighborhood of $\bar{u}$. Hence, multiplying relation \eqref{srd} by minus one, we may apply the same Clarke's theorem in order to derive the inclusion
\begin{equation}\label{clarkeinclu2}
\partial^C(-\varphi )(\bar{u})\subseteq\int\limits_{\mathbb{S}^{m-1}}\partial^C_u(-e)(\bar{u},v) \D \mu_\zeta (v).
\end{equation}
Now, the addendum to Clarke's theorem mentioned above yields that equality in \eqref{clarkeinclu2} and Clarke regularity of $-\varphi$ would hold true if the functions $-e(\cdot ,v)$ were Clarke regular at $\bar{u}$. Contrary to the functions $e(\cdot ,v)$ themselves, this will hold true for $-e(\cdot ,v)$, so that we can conclude equality instead of inclusion in \eqref{clarkeinclu} upon multiplying the relation \eqref{clarkeinclu2} with minus one. 
Summarizing, the assertion of Theorem \ref{clarkeequality} will hold true as a consequence of the following 
\begin{lemma}
Under the assumptions of Theorem \ref{clarkeequality}, for each $v\in\mathbb{S}^{m-1}$, the functions $-e(\cdot ,v)$ are Clarke regular at $\bar{u}$.
\end{lemma}
\proof
We define functions $e_w\colon U\times\mathbb{S}^{m-1}\to\mathbb{R}$ for $w\in K$ by $e_w(u,v):=\mu_\chi (A(u,v,w))$, where $A(u,v,w):=\{r\geq 0\mid h(u,r\Sigma^{1/2}v,w)\leq 0\}$. By continuity of $g$ and $g(\bar{u},0)<0$, there exists a neighborhood $\mathcal{N}$ of $\bar{u}$ such that $g(u,0)<0$ or $h(u,0,w)<0$ for all $u\in\mathcal{N}$ and all $w\in K$. Then, the convexity of $h$ in the second argument yields  the existence of a positive (possibly infinite) number $\rho (u,v,w)$ such that $A(u,v,w)=[0,\rho (u,v,w)]$ for all $(u,v,w)\in\mathcal{N}\times\mathbb{S}^{m-1}\times K$. Now, the definitions of $e$ and $g$ provide that, for all $(u,v)\in\mathcal{N}\times\mathbb{S}^{m-1}$, 
\begin{eqnarray*}
e(u,v)&=&\mu_\chi\left(\bigcap_{w\in K}A(u,v,w)\right)=\mu_\chi \left(\bigcap_{w\in K}[0,\rho (u,v,w)]\right)\\&=&
\mu_\chi\left([0,\inf\limits_{w\in K}\rho (u,v,w)]\right)=\inf\limits_{w\in K}\mu_\chi ([0,\rho (u,v,w)])=\inf\limits_{w\in K}e_w(u,v). 
\end{eqnarray*}
Here, the penultimate identity follows from $\mu_\chi$ being an absolutely continuous measure.
Condition 4. of Theorem \ref{clarkeequality} ensures that $e_w(\cdot ,v)$ is continuously differentiable on $\mathcal{N}$ for each $w\in K$ and each $v\in\mathbb{S}^{m-1}$ by \cite[Corollary 3.2 and Example 3.1]{ackooij23} and its  derivative is given by (with $\rho$ introduced above)
\begin{align}\label{gra_form}
    D_u e_w(u,v) =\left\{  \begin{array}{cc}
   -f_\chi (\rho (u,v,w)) \frac{ D_u h(u,\rho (u,v,w) \Sigma^{1/2}v,w) }{ \langle D_z h(u,\rho (u,v,w) \Sigma^{1/2}v,w), \Sigma^{1/2}v \rangle  }      &  \text{ if } \rho (u,v,w)< \infty \\
     0     & \text{ if } \rho (u,v,w)= \infty .
    \end{array} \right.
\end{align}
Here, $f_\chi$ denotes the density of the Chi- distribution $\mu_\chi$ with $m$ degrees of freedom. The Clarke regularity of $-e(\cdot ,v)=\sup_{w\in K}-e_w(u,v)$ at $\bar{u}$ for every $v \in \mathbb{S}^{m-1}$ can be checked by applying \cite[Theorem 2.8.2]{clarke1983}. Indeed, according to this theorem and translated to our setting, it will be sufficient to show the following statements for each $v\in \mathbb{S}^{m-1}$ separately:
\begin{enumerate}[label=(\alph*)]
\item The mapping $\mathcal{N} \times K \ni (u,w) \mapsto D_u (-e_w)(u,v)$ is continuous
\item The functions $e_{w}(\cdot,v)$ with $w\in K$ are Lipschitz continuous on a neighborhood of $\bar{u}$ with some common modulus;
\item The set $\{e_w(\bar{u},v)\mid w\in K\}$ is bounded;
\item The mapping $w\mapsto e_w(u,v)=\mu_\chi ([0,\rho (u,v,w)])$ is continuous at each $u\in\mathcal{N}$.
\end{enumerate}
In order to show (a), consider a sequence $(u_k,w_k) \to (u,w)$. On the one hand, if $\rho (u,v,w)<\infty$, we have that the continuity of $\rho$ \cite[Lemma 4.10 (ii)]{MR4000225} and the continuity of $D_uh, D_zh$ by condition 2. of Theorem \ref{clarkeequality}~imply via \eqref{gra_form} that $D_u (-e_{w_k})(u_k,v) \to  D_u (-e_{w})(u,v) $. On the other hand, if $\rho(u,v,w)=\infty$, then again by \cite[Lemma 4.10 (ii)]
{MR4000225}, we have that $\rho (u_k,v,w_k)\to\infty$. Now, adapting \cite[Lemma 3.2]{ackooij23} to our setting, we infer for $k$ large enough the estimate
\[
\|D_u e_{w_k}(u_k,v)\| \leq\frac{c }{|h(u,0,w_k)|} f_\chi (\rho (u_k,v,w_k))  \rho (u_k,v,w_k) \exp\left(   4\rho (u_k,v,w_k) \|\Sigma^{1/2}\|\right),
\]
where $c$ is the constant from condition 4. of Theorem \ref{clarkeequality}. For $k\to\infty$, the expression $|h(u,0,w_k)|$ tends to $|h(u,0,w)|$, which is nonzero by definition of $\mathcal{N}$ (see above). Moreover, given the explicit formula for the density of the Chi-distribution with $m$ degrees of freedom (including some normalizing constant $\tilde{c}>0$), we get that
\[
f_\chi (t)t\exp\left(4t\|\Sigma^{1/2}\|\right)=\tilde{c}t^m
\exp(-t^2/2)\exp\left(4t\|\Sigma^{1/2}\|\right)\to 0 \quad (\textup{as } t\to \infty).
\]
Hence, by \eqref{gra_form}, $D_u e_{w}(u,v)=0=\lim_{k\to\infty}D_u e_{w_k}(u_k,v)$, which shows the continuity of $(u,w) \mapsto D_u (-e_w)(u,v)$. To proceed with items (b-d), observe that, from the already shown continuity of the mapping $(u,w) \mapsto D_u e_{w}(u,v)$, it follows that the function $u\mapsto\max_{w\in K}\|D_u e_{w}(u,v)\|$ is continuous on $\mathcal{N}$. Hence, there is a neighborhood 
$\mathcal{N}'$ of $\bar{u}$ and some $\kappa >0$ such that
$\|D_u e_{w}(u,v)\|\leq\kappa$ for all $w\in K$ and all $u\in\mathcal{N}'$. This  implies (b) with the common Lipschitz modulus $\kappa$, while (c) is trivial because $e_w$ is a probability. Concerning (d), fix an arbitrary $(u,v,w)\in\mathcal{N}\times\mathbb{S}^{m-1}\times K$ and a sequence $w_k\to w$ with $w_k\in K$. We repeat the case distinction on $\rho (u,v,w)$ from above: If $\rho (u,v,w)<\infty$, then $\rho (u,v,w_k)\to\rho (u,v,w)$ by the continuity of $\rho$ already mentioned above. Moreover, $\rho (u,v,w_k)<\infty$ for $k$ large enough and, hence, with $F_\chi$ denoting the continuous distribution function of the Chi-distribution with $m$ degrees of freedom, we arrive at
\begin{eqnarray*}
e_{w_k}(u,v)&=&\mu_\chi ([0,\rho (u,v,w_k)])=F_\chi (\rho (u,v,w_k))\to_k F_\chi (\rho (u,v,w))\\&=&\mu_\chi ([0,\rho (u,v,w)])=e_w(u,v),
\end{eqnarray*}
which is the desired continuity in the first case. Otherwise, if $\rho(u,v,w)=\infty$, then, as already seen above, we have that $\rho (u,v,w_k)\to\infty$. Let $\varepsilon>0$ be arbitrary.
Since $\lim_{t\to\infty}F_\chi (t)=1$, it follows that $e_{w_k}(u,v)=F_\chi (\rho (u,v,w_k))\geq 1-\varepsilon$ whenever $k$ is large enough and $\rho (u,v,w_k)<\infty$. If, in contrast, $\rho (u,v,w_k)=\infty$ for some $k$, then 
\[
e_{w_k}(u,v)=\mu_\chi ([0,\infty ))=1.
\]
Hence, $e_{w_k}(u,v)\geq 1-\varepsilon$ whenever $k$ is large enough. Consequently,
\[
e_{w_k}(u,v)\to_k 1=\mu_\chi ([0,\infty ))=\mu_\chi ([0,\rho(u,v,w) )=e_w(u,v).
\]
This proves the desired continuity in the second case and, hence, (d).

\qed
\bibliography{references}
\end{document}